\documentclass{IEEEtran}
\usepackage{textcomp}
\usepackage{amsmath,amsfonts,amssymb,amscd}
\usepackage[colorlinks,linkcolor=blue,citecolor=blue]{hyperref}
\usepackage[font=small]{caption}
\usepackage{latexsym}
\usepackage{graphicx,subfigure}
\usepackage{hyperref}
\usepackage{color}
\usepackage{bm}
\usepackage{cite}
\usepackage{algorithm}
\usepackage{algorithmic}
\usepackage{mathrsfs}
\usepackage{subeqnarray}
\usepackage{cases}
\usepackage{arydshln}
\usepackage{multirow}
\usepackage{geometry}
\usepackage{extarrows}

\def\BibTeX{{\rm B\kern-.05em{\sc i\kern-.025em b}\kern-.08em
		T\kern-.1667em\lower.7ex\hbox{E}\kern-.125emX}}
\newtheorem{thm}{Theorem}
\newtheorem{lem}{Lemma}
\newtheorem{remark}{Remark}
\newtheorem{ass}{Assumption}
\newtheorem{definition}{Definition}
\newtheorem{prop}{Proposition}

\newtheorem{corollary}{Corollary}

\newcommand{\norm}[1]{\left\lVert#1\right\rVert}

\begin{document}

\title{Fixed-Time Gradient Flows for Solving Constrained Optimization: A Unified Approach}

\author{Xinli~Shi,~\IEEEmembership{Senior Member, IEEE}, Xiangping Xu, Guanghui~Wen,~\IEEEmembership{Senior Member, IEEE}, Jinde~Cao,~\IEEEmembership{Fellow,~IEEE}
\thanks{Citation: X. Shi, X. Xu, G. Wen, and, J. Cao, “Fixed-time gradient flows for solving constrained optimization: A unified approach,” IEEE/CAA Journal of Automatica Sinica, doi: 10.1109/JAS.2023.124089, 2023. }
\thanks{X.~Shi is with the School of Cyber Science and Engineering, Southeast University, Nanjing 210096, China (e-mail: xinli\_shi@seu.edu.cn).}
\thanks{X.~Xu, G. Wen and J. Cao are with the School of Mathematics, Southeast University, Nanjing 210096, China (e-mail: xpxu2021@seu.edu.cn; ghwen@seu.edu.cn; jdcao@seu.edu.cn).}
}




\maketitle

\begin{abstract}
The accelerated	method in solving optimization problems has always been an absorbing topic. Based on the fixed-time (FxT) stability of nonlinear dynamical systems, we provide a unified approach for designing FxT gradient flows (FxTGFs). First, a general class of nonlinear functions in designing FxTGFs is provided. A unified method for designing first-order FxTGFs is shown under Polyak-\L jasiewicz inequality assumption, a weaker condition than strong convexity. When there exist both bounded and vanishing disturbances in the gradient flow, a specific class of nonsmooth robust FxTGFs with disturbance rejection is presented. Under the strict convexity assumption, Newton-based FxTGFs is given and further extended to solve time-varying optimization. Besides, the proposed FxTGFs are further used for solving equation-constrained optimization. Moreover, an FxT proximal gradient flow with a wide range of parameters is provided for solving nonsmooth composite optimization. To show the effectiveness of various FxTGFs, the static regret analysis for several typical FxTGFs are also provided in detail. Finally, the proposed FxTGFs are applied to solve two network problems, i.e., the network consensus problem and solving a system linear equations, respectively, from the respective of optimization. Particularly, by choosing component-wisely sign-preserving functions, these problems can be solved in a distributed way, which extends the existing results. The accelerated convergence and robustness of the proposed FxTGFs are validated in several numerical examples stemming from practical applications.
\end{abstract}


\IEEEpeerreviewmaketitle

\section{Introduction}
The gradient dynamical systems in virtue of ordinary differential equations (ODE) for solving unconstrained optimization have gained a lot of interest. Generally, compared with the discrete-time scheme, the continuous-time methods often have better intuition and simpler structures. Besides, based on the Lyapunov method from control theory, the stability of equilibria referring to the optimal solutions can be shown directly. For example, given a locally convex function $f: \mathbb{R}^n\rightarrow \mathbb{R}$, one can prove that the classic gradient flow (GF): $\dot{x}=-\nabla f(x)$ will converge to the strict minima of $f$ asymptotically if the level sets of $f$ are compact. GF and its variants have been applied to many problems including distributed optimization \cite{Wang2011}\cite{Liu2017}, network consensus \cite{cortes2006} and image processing \cite{Xu1998}. 

Nowadays, accelerated methods for solving optimization have aroused more attention to obtain a fast convergence rate. In \cite{Su2014}, a second-order ODE is derived to model Nesterov’s accelerated
algorithms and proven to have a linear convergence rate when the cost function is strongly convex. In \cite{Wibisono2016}, Bregman Lagrangian is introduced to generate a class of accelerated dynamics with exponential convergence rates. More second-order ODE methods for solving unconstrained optimization can be found in \cite{Attouch2018,Vassilis2018,Sebbouh2020,Serna2021}. Although these accelerated methods can achieve better convergence rates than the classic GF, the optimal solution will still be reached in infinite time despite the polynomial or exponential convergence rate, which limits the applications in scenarios where the optimization is required to be solved in finite or prescribed time.     

To achieve finite-time (FT) convergence, several nonlinear GFs have been proposed by showing the FT stability of the provided dynamics based on the seminal work \cite{Bhat2000}. FT stability can guarantee that the dynamical system converges to the equilibrium in finite time that may depend upon the initial states. FT control of nonlinear systems has been investigated in \cite{Yu2005,Shen2012, Aouiti2021FT}. For solving optimization, two discontinuous FT gradient flows (FTGFs) are studied in \cite{cortes2006} with application in the network consensus problem. The work \cite{Romero2020} studies two FTGFs, i.e., re-scaled gradient flow and signed-gradient flow, which extend the discontinuous protocols \cite{cortes2006} to continuous versions. The work \cite{Chen2020sign} proposes a sign projected FTGF for solving linear equation-constrained optimization. In \cite{Wei2022}, fractional-order gradient dynamics are proposed with FT convergence. Besides centralized optimization, FT algorithms also contribute to fast distributed optimization. For example, several works combine the FTGFs with consensus protocols for solving distributed optimization \cite{Zhou2019tac,Shi2022cyber,Shi2022TAC,Shi2023}. In \cite{Zhou2019tac}, a nonsmooth projection-based distributed protocol is provided for solving a system of linear equations with minimum $l_1$-norm in finite time. Recently, a discontinuous primal-dual dynamics with FT convergence has been proposed in \cite{Shi2022cyber} for solving constrained (distributed) optimization. In \cite{Shi2022TAC}, several discontinuous continuous-time dynamics with bounded inputs are proposed to solve constrained distributed optimization based on FT average tracking dynamics.

As the finite settling time depends upon the initial condition, it can be unbounded as the initial state increases. When the bound of the finite setting time is obtained regardless of the initial conditions, the fixed-time (FxT) stability is established \cite{Polyakov2012}. A series of FxT control methods have been surveyed in \cite{LiuCAA2022}, and FxT stabilization/synchronization of neural networks has been investigated in \cite{Aouiti2019,Alimi2019,Cao2017, Aouiti2020,Aouiti2020CSSP,Aouiti2021}. For solving optimization problems, in \cite{Garg2021tac}, a specific class of FxT gradient flows (FxTGFs) with continuous right-hand sides is proposed. Furthermore, in \cite{Budhraja2021}, the robustness of the FxTGFs given in \cite{Garg2021tac} is analyzed with vanishing disturbance and the regret analysis is further provided under the Polyak-\L jasiewicz (PL) inequality assumption. However, the more general disturbance is not considered and a detailed analysis of the regret bound is not given. As for the composite optimization, an FxT stable proximal dynamics is provided in \cite{Garg2019MVI}, where the smooth term is required to be strongly convex and the related parameters are restrained. Moreover, in \cite{Ju2021Tcyber}, a novel FxT convergent neurodynamics is given with the strong pesudomonotone operator. For the linear equation problem with $L_1$-minimization, the work \cite{HeTNNLS2021} provides a modified projection neural network with fixed-time convergence. As for the network optimization, in \cite{Wu2022SMC}, an FxT convergent algorithm is given for solving distributed optimization based on a zero-gradient-sum method. In \cite{Shi2020auto}, several FT/FxT nonlinear algorithms are designed for solving linear equations based on projection methods.


In this paper, we provide a unified approach for designing FxTGFs under PL assumption, which generalizes the existing FT/FxT methods \cite{cortes2006,Garg2021tac,Guo2022CAA,Budhraja2021,Garg2019MVI}. PL inequality is more general than strong convexity and has been reported as the weakest condition used for achieving linear convergence \cite{Karimi2016}\cite{Yi2022CAA}. The provided scheme not only covers the existing algorithms but also can generate some novel FxTGFs with discontinuous or continuous right-hand sides. We further investigate the robustness of a special class of FxTGFs in the presence of both bounded and vanishing disturbances, and the FxTGFs for constrained optimization. The contributions of this work compared with previous methods can be summarized as follows.
\begin{enumerate}
\item We first provide a general class of nonlinear protocols for designing the first-order FxTGFs from a unified framework based on nonsmooth analysis and PL inequality. Moreover, with both bounded and vanishing disturbances, a specific class of nonsmooth robust FxTGFs (RFxTGFs) with disturbance rejection is given, extending the work \cite{Budhraja2021} where only vanishing disturbance is considered. Besides, a class of Newton-based FxTGFs (NFxTGFs) is provided and further extended to solve time-varying optimization. The fixed settling time bounds are further given explicitly for all FxTGFs. 
	
\item We provide a projected FxTGF (PFxTGF) for solving equation-constrained optimization without the full row-rank assumption, differently from the work \cite{Garg2021tac} where the dual and primal FxTGFs are used sequentially requiring more computation complexity. In addition, an FxT 
proximal gradient flow (FxTPGF) is also provided for solving nonsmooth composite optimization based on proximal PL inequality and forward–backward (FB) envelope, relaxing the parameter conditions in \cite{Garg2019MVI} considerably.

\item To investigate the performance of various FxTGFs, estimations of the static regret bounds for several typical FxTGFs are provided with a more detailed analysis compared with \cite{Budhraja2021}, and we conclude that not all FxTGFs have fixed regret bounds regardless of the initial condition.
\item As applications, the proposed FxTGFs are used to solve two network problems, i.e., network consensus problem and solving a system of linear equations, extending the existing results \cite{Xiao2009, Shi2020TAC, Shi2019TNSE,Zuo2014}. Particularly, by choosing component-wisely sign-preserving functions, these problems can be solved in a distributed way. 
\end{enumerate}

\begin{figure}[t]
	\centering
	{\includegraphics[width=.45\textwidth]{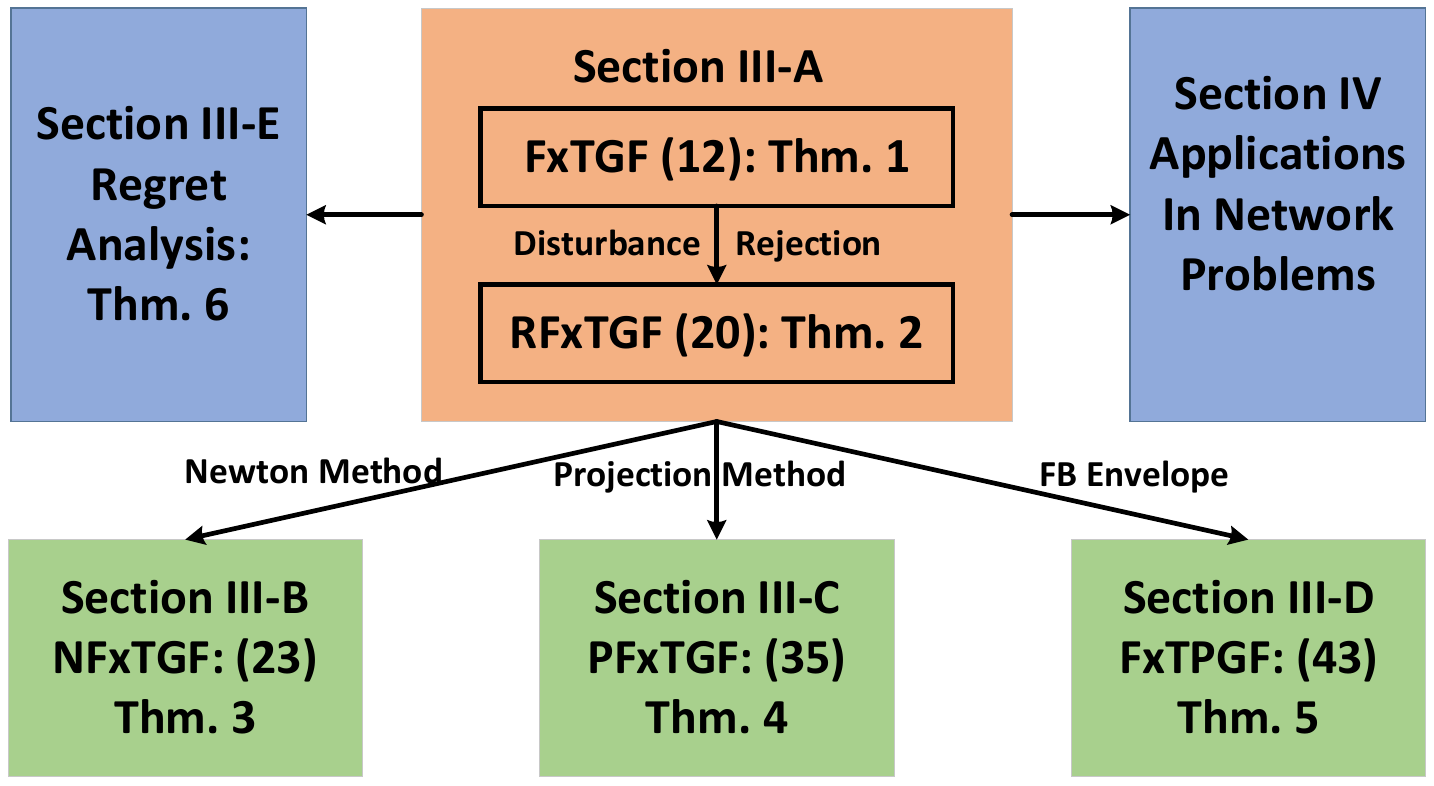}}
	\caption{The overview of the proposed approaches.}
	\label{overview}
\end{figure}

The organization of this work is given as follows. In Section \ref{Preliminaries}, the preliminary notations and concepts are introduced, as well as the problem formulation. The overview of the proposed methods is shown in Fig. \ref{overview}. Section \ref{FxTGFs} presents the main results including the first-order FxTGFs and NFxTGFs design, RFxTGFs with disturbance rejection, PFxTGFs for equation-constrained optimization, FxTPGFs for composite optimization and regret analysis. Moreover, the proposed FxTGFs are applied to two network problems in Section \ref{application}. In Section \ref{Numberical}, several practical case studies are conducted for testifying the proposed protocols with and without disturbances. Finally, conclusions are drawn in Section \ref{conclusion}.

\par
\section{Preliminaries and Problem Formulation}\label{Preliminaries}

\subsection{Notation}
The set of $n$-dimensional vectors is represented by $\mathbb{R}^n$. We denote $\bm{1}_n$ $\in \mathbb{R}^n$ as the vector with elements being all one, and $I_n$ as an $n$-dimensional identity matrix. Let $\langle n\rangle=\{1,2,\cdots,n\}$. For $x=[x_{1},\ldots,x_{n}]^{T}\in \mathbb{R}^n$ and any $p>0$, the $p$-norm of $x$ is $\norm{x}_p=(\sum_{i=1}^n |x_i|^p)^{\frac{1}{p}}$, and we define $|x|^p = [|x_i|^p]_{i\in \langle n \rangle}$, $e^{|x|} = [e^{|x_i|}]_{i\in \langle n \rangle}$ and $\text{sign}(x)=[\text{sign}(x_1),...,
\text{sign}(x_n)]^T$ with $\text{sign}(\cdot)$ being the traditional signum function.
For any pair of vectors $x,y \in \mathbb{R}^n$, we define the operator $x\odot y = [x_iy_i]_{i\in \langle n\rangle}$. 

For a matrix $A \in \mathbb{R}^{n\times n}$, the null and image space of the linear map $A$ is respectively represented by $\mathcal{N}(A)$ and $\mathcal{R}(A)$. When $A$ is positive semi-definite, $\lambda_2(A)$ denotes the smallest nonzero eigenvalue of $A$. A differentiable function $f(x): \mathbb{R}^n\rightarrow \mathbb{R}$ is $\mu-$strongly convex if it satisfies $f(y)\geq f(x)+\nabla^T f(x)(y-x)+ \frac{\mu}{2}\|x-y\|_2^2$ for any $x,y \in \mathbb{R}^n$ and some $\mu>0$. A vector-valued function $f:\mathbb{R}^n\rightarrow \mathbb{R}^n$ is called \textit{sign-preserving} if $\text{sign}(f(x))=\text{sign}(x), \forall x\in \mathbb{R}^n$. If the sign-preserving function $f(x)$ is component-wisely defined, i.e., $f(x)=[f_i(x_i)]_{i\in \langle n \rangle}$, it is called \textit{component-wisely sign-preserving} \cite{Shi2022cyber}.

\subsection{Nonsmooth Analysis}\label{Nonsmooth}
Consider the following differential system
\begin{align}\label{auto_dynamic}
\dot{x}(t)=X(x(t)), x(t_0)=x_0
\end{align}
with $x\in \mathbb{R}^n$. When the map $X: \mathbb{R}^n \rightarrow \mathbb{R}^n$ is discontinuous, the Filippov solution to (\ref{auto_dynamic}) is investigated, which is an absolutely continuous map $\phi:I\subset \mathbb{R} \rightarrow \mathbb{R}^n $ satisfying 
\begin{align}\label{filippov}
\dot{\phi}(t) \in \mathcal{F}[X](\phi(t)).
\end{align}
The definition of Filippov set-valued map $\mathcal{F}[X](\cdot)$ is given by 
\begin{align}\label{setmap}
\mathcal{F}[X](x) \triangleq  \bigcap_{\delta>0} \bigcap_{\omega(\mathcal{N}_0)=0}\overline{co}\{X(x^{\delta} \backslash \mathcal{N}_0)\}, 
\end{align}
where $x^{\delta}$ is an open ball centered at $x$ with radius $\delta>0$, and $\omega(\mathcal{N}_0)$ denotes the Lebesgue measure of $\mathcal{N}_0$. 
As known, when $f$ is measurable and locally essentially bounded (abbreviated as m.l.b.), a Filippov solution to \eqref{auto_dynamic} always exists \cite{Cortes2008}. 

For the sake of nonsmooth analysis, we introduce the following concepts. The \textit{set-valued Lie derivative} $\widetilde{\mathcal{L}}_{\mathcal{F}}V$ of a locally Lipschitz continuous map $V: \mathbb{R}^n \rightarrow \mathbb{R}$ related to a set-valued map $\mathcal{F}$ at $x$ is defined as
\begin{align*}
	\widetilde{\mathcal{L}}_{\mathcal{F}}V(x) \triangleq \{a \in \mathbb{R}| \exists \eta \in \mathcal{F}(x) \ \text{with} \ \nu^T \eta=a, \forall \nu \in  \partial V(x)\}
\end{align*}
with $\partial V(x)$ being the generalized gradient of $V$ \cite{Cortes2008}.
The following lemma gives the possible dynamics of $V$ along the solutions to \eqref{filippov}.
\begin{lem}\cite[Prop. 10]{Cortes2008}\label{lem-non}
For a locally Lipschitz continuous and regular map $V: \mathbb{R}^n \rightarrow \mathbb{R}$ and a Filippov solution $\phi(t): I\subset \mathbb{R} \rightarrow \mathbb{R}^n$ to differential inclusion \eqref{filippov}, $V(\phi(t))$ is absolutely continuous and satisfies 
\begin{align}
  \frac{dV(\phi(t))}{dt} \overset{a.e.}{\in} \widetilde{\mathcal{L}}_{\mathcal{F}[X]}V(\phi(t)), \ t\in I.
\end{align}
\end{lem}
\begin{lem}\label{lem-fixed}
Consider the system (\ref{auto_dynamic}) with $X(0)=0$.	Let $x(t; x_0)$ be any solution of the system (\ref{auto_dynamic}) starting from $x_0\in  \mathbb{R}^n$. Suppose that there exists a continuous radially unbounded and positive definite function $V(x)$ satisfying $V(x)=0 \Rightarrow x=0$. 
	\begin{enumerate}
		\item \cite[Theorem 4.2]{Bhat2000} If $\dot{V}(x(t))\leq -cV^r(x(t))$ with $c>0, r\in (0,1)$, the origin is FT stable with settling time $T(x_0)\leq \frac{V^{1-r}(x_0)}{c(1-r)}$.
		\item \cite[Lemma 1]{Polyakov2012} If $V(x)$ satisfies that
		\begin{align}
			\dot{V}(x(t))\leq -(c_1 V^{r_1}(x(t))+ c_1 V^{r_2}(x(t) )^k
		\end{align}
		for some positive scalars $c_1,c_2,r_1,r_2: 0<r_1k<1, r_2k>1$, the origin is FxT stable with settling time estimated by
		\begin{align}
			T(x_0)\leq \frac{1}{c_1^k(1-r_1k)} + \frac{1}{c_2^k(r_2k-1)}.
		\end{align}
	\end{enumerate}
\end{lem}

\begin{ass}\label{ass-f}
	There exists an optimizer $x^*$ such that the cost function $f(x)$ attains its minimum value $f^*=f(x^*)>-\infty$. 
\end{ass}
\begin{ass}\label{ass-f2}
	The cost function $f(x)$ is continuously differentiable and satisfies Polyak-\L jasiewicz (PL) inequality with the parameter $\mu>0$, i.e.,
	\begin{align}\label{PL}
		\frac{1}{2}\|\nabla f(x)\|_2^2 \geq \mu (f(x)-f^*), \forall x\in \mathbb{R}^n.
	\end{align} 
\end{ass}

\subsection{Problem Formulation}\label{Notations}
Consider the following optimization 
\begin{align}\label{CO}
	\min_{x\in \mathcal{X}} f(x),
\end{align}
where $f:\mathbb{R}^n\rightarrow \mathbb{R}$ is the cost function to be minimized and $\mathcal{X}\subseteq \mathbb{R}^n$ is a convex closed set. Denote the set of optimal solutions to \eqref{CO} by $\mathbb{X}^*$, which is supposed to be nonempty. We further make Assumptions \ref{ass-f} and \ref{ass-f2} for \eqref{CO}. PL inequality is the weakest condition used in literature for achieving linear or exponential convergence \cite{Karimi2016}. Note that any strongly convex function satisfies Assumption \ref{ass-f2}. Besides, when the matrix $P$ is positive semi-definite, $f(x)=x^TPx$ satisfies Assumption \ref{ass-f2}. According to \cite{Karimi2016}, PL inequality \eqref{PL} implies that $f(x)$ has a quadratic growth, i.e.,
\begin{align}\label{QG}
			 f(x)-f^*\geq 2\mu \|x-[x]^*\|_2^2, \forall x\in \mathbb{R}^n
\end{align} 
where $[x]^*$ denotes the projection of $x$ on $\mathbb{X}^*$. Combining with \eqref{PL}, we derive that 
\begin{align}\label{QG2}
\|\nabla f(x)\|_2 \geq 2\mu \|x-[x]^*\|_2, \forall x\in \mathbb{R}^n.
\end{align} 

In this note, we aim to design the FxT convergent gradient protocol $u: \mathbb{R}^n\rightarrow \mathbb{R}^n$ with disturbance rejection (when disturbance exists) such that the following dynamics
\begin{align}\label{pro-u}
	\dot{x} = -u(x) +d(x,t)
\end{align}
can reach $\mathbb{X}^*$ in a fixed time regardless of the initial states. For convenience, we denote the right-side of \eqref{pro-u} by $X_{u,d}: \mathbb{R}^{n+1} \rightarrow \mathbb{R}^{n}$. Specifically, when $d(x,t)\equiv 0$, the dynamics \eqref{pro-u} is referred as a nominal system.

\section{FxTGFs: Methodology, Disturbance Rejection and Regret Analysis}\label{FxTGFs}
In this section, we first give a unified scheme for designing first-order FxTGFs \eqref{pro-u} without and with disturbances, respectively. A class of second-order FxTGFs is also provided. Then, several FxTGFs are presented for solving equation-constrained optimization, and an FxTPGF is presented for solving nonsmooth composite optimization. Finally, the regret analyses of several typical FxTGFs for unconstrained optimization are given to measure the performance of various continuous-time algorithms. Several frequently used notations and abbreviations in this part are summarized in Table \ref{tb-notations}. 

\subsection{First-Order FxTGFs For Unconstrained Optimization}\label{subsection-gn}
At the beginning, we will give a unified algorithm design for FxTGFs to solve \eqref{CO} in the absent of disturbances, i.e., the nominal system \eqref{pro-u} with $d(x,t)\equiv 0$. As preparation, we first give Definition \ref{g-def}. To achieve FxT convergence, we choose
\begin{align}\label{prot-gpq}
	u(x(t))= g(\nabla f(x)) = g_p(\nabla f(x))+g_q(\nabla f(x)), 
\end{align}
which satisfies Assumption \ref{ass-gn}. Note that when $g$ is continuous at $y \in \mathbb{R}^n$, $\mathcal{F}[g](y) =\{g(y)\}$. In Assumption \ref{ass-gn}, the subfunction $g_p$ is used for ensuring FT stability, and $g_q$ is additionally required for FxT stability.

\begin{table}[t]
	\caption{Typical notations and abbreviations.}
	\label{tb-notations}
	\centering
	\small
	\begin{tabular}{p{1.8cm}p{1.4cm}p{4.7cm}}
		\hline
		Notation& Location     &  Description  \\
		\hline
	    $\mathbb{X}^*$& -  & set of optimal solutions to \eqref{CO} \\  
	    $[x]^*$& -  & projection of $x$ on $\mathbb{X}^*$ \\  
	     $f^*$& -  & optimal objective value to \eqref{CO} \\  
		$g=g_p+g_q$& Ass. \ref{ass-gn}   & FxT stable function used in FxTGFs \\  
		$d(x,t)$ &  Ass. \ref{ass-d}   &   system disturbance \\ 
		$ \mathcal{F}[g](y)$& Eq. \eqref{setmap} &  Filippov set-valued map of $g$ \\
		$\mathcal{F}_p^{\sigma}$& Def. \ref{g-def} & set of functions satisfying \eqref{gp} \\
		$\partial \|y\|_r $& Eq. \eqref{grad-norm0} & (sub)gradient of norm $\|y\|_r, r\geq 1$\\
		$T(x_0)$& Theorems  &   FT/FxT settling time from $x_0$ \\ 
		$\mathcal{R}/\mathcal{N}(A)$ &  Sec. \ref{sub-COLE}   & image/null space of linear map $A$\\
		$P/P_A$ & Sec. \ref{sub-COLE}          & projection matrix onto $\mathcal{N}(A)$\\
		$\lambda_2(AA^T)$ & Sec. \ref{sub-COLE}          & smallest nonzero eigenvalue of $AA^T$\\
		$\text{prox}_{\lambda h}(x)$ &  Eq. \eqref{prox}  & proximal operator of $h$ with $\lambda>0$ \\
		$\mathcal{R}(T,x_0)$&  \eqref{regret-fun} & regret function over $[0 \ T]$ from $x_0$ \\
		\hline  
		PL inequality  & \eqref{PL} & General assumption for objective\\
		FxTGF  & \eqref{prot-gpq} & First-order FxTGF\\
		RFxTGF &  \eqref{prot-gd} & Robust FxTGF\\
		NFxTGF  & \eqref{prot-gn2} & Newton-based FxTGF\\
	    PFxTGF  & \eqref{prot-free} & Projected FxTGF \\
	    FxTPGF  & \eqref{u-cop1} & FxT Proximal GF \\
		\hline
	\end{tabular}
\end{table}

\begin{definition}\label{g-def}
For a measurable, locally essentially bounded and sign preserving function $g_p: \mathbb{R}^n\rightarrow \mathbb{R}^n$ with $p \geq 0$, we denote $g_p \in \mathcal{F}_p^{\sigma} $ if there exists $\sigma>0$ such that it satisfies 
\begin{align}\label{gp}
	\min_{\eta \in \mathcal{F}[g_p](y)} \eta^Ty\geq \sigma \|y\|_2^{1+p},  \forall y \in \mathbb{R}^n.
\end{align}
\end{definition}

\begin{ass}\label{ass-gn}
	The function $g=g_p+g_q: \mathbb{R}^n\rightarrow \mathbb{R}^n$, where $g_p \in \mathcal{F}_p^{\sigma}$ and $g_q \in \mathcal{F}_q^{\rho}$ with $p \in [0,1), q>1$ and $\sigma, \rho>0$.
\end{ass}

\begin{lem}\cite[Thm. 1]{Goldberg1987}\label{lem-ineq}
	For any $x\in \mathbb{R}^n$ and $r>s\geq1$, it holds that
	\begin{align}
		\|x\|_r\leq \|x\|_s \leq n^{\frac{1}{s}-\frac{1}{r}}\|x\|_r.
	\end{align}
\end{lem}

Based on Lemma \ref{lem-ineq}, one can obtain several typical functions/subfunctions that satisfy Assumption \ref{ass-gn}, as presented in Lemma \ref{lem-g}.
\begin{lem}\label{lem-g}
	For $p\in [0,1)$, $q>1$ and $r\in [1, +\infty]$, it holds that
	\begin{enumerate}
		\item $g_0(y) = \partial \|y\|_r \in \mathcal{F}_0^{\sigma}$ with $\sigma=1, r \in [1,2]$ and $\sigma=n^{\frac{1}{2}-\frac{1}{r}}$ when $r>2$;
		\item $g_p(y) =  \frac{y}{\|y\|_r^{1-p}} \in \mathcal{F}_p^{\sigma}$ with $\sigma=n^{\frac{1-p}{2}-\frac{1-p}{r}}, r \in [1,2]$ and $\sigma=1$ when $r>2$;
		\item $g_q(y)=y\|y\|_r^{q-1} \in \mathcal{F}_q^{\rho}$ with $\rho=1, r\in [1,2]$ and $\rho=n^{(q-1)(\frac{1}{2}-\frac{1}{r})}$ when $r>2$;
		\item $g_p(y)=\text{sign}(y)\odot |y|^{p}\in \mathcal{F}_p^{1}$ and $g_q(y)=\text{sign}(y)\odot |y|^{q}\in \mathcal{F}_q^{\rho}$ with $\rho=n^{(1-\frac{q+1}{2})}$;
	    \item $g_{e,2}(y)=\frac{ye^{\|y\|_2}}{\|y\|_2}$ and $g_{e,1}(y)=\text{sign}(y)\odot e^{|y|}$ satisfy Assumption \ref{ass-gn}. Moreover, for any positive integer $k$,
	    \begin{align}\label{ge-property}
	    	y^Tg_{e,2}(y) =\frac{y^Tye^{\|y\|_2}}{\|y\|_2}=\|y\|_2e^{\|y\|_2} \geq  \frac{\|y\|_2^{k+1}}{k!},
	    \end{align}
    and $y^Tg_{e,1} \geq (\frac{y}{\sqrt{n}})^Tg_{e,2}(\frac{y}{\sqrt{n}})$, i.e.,
    \begin{align}\label{ge-property2}
        \sum_{i=1}^{n} |y_i|e^{|y_i|} \geq \frac{\|y\|_2}{\sqrt{n}}e^{\frac{\|y\|_2}{\sqrt{n}}}
    \end{align} 
	\end{enumerate}
\end{lem}
\begin{proof}
See Appendix A.
\end{proof}

More functions satisfying Assumption \ref{ass-gn} are summarized in Table \ref{tb-g}, where $\partial \|y\|_r$ represents the (sub)gradient of $\|y\|_r, r\geq 1$ w.r.t. $y$, i.e., 
\begin{align}\label{grad-norm0}
\partial \|y\|_r =   \left\{\begin{aligned}
		\frac{\text{sign}(y)\odot|y|^{r-1}}{\|y\|_r^{r-1}}, &\  y\neq 0, \\
		0 \qquad \qquad, & \ y=0.
	\end{aligned}\right.
\end{align}
Besides, as the function $g(y)=\text{sign}(y)\odot |y|^{\alpha}$ or $\text{sign}(y)\odot e^{|y|}$ is defined component-wisely, they can be used for designing distributed network protocols. By Assumptions \ref{ass-f2} and \ref{ass-gn}, $u(x(t))$ is locally essentially bounded and measurable, which implies that the  Filippov solution to the nominal system \eqref{pro-u} always exists.
With the previous preparations, we give the following result for the nominal system \eqref{pro-u}.

\begin{table}[t]
	\caption{Typical functions satisfying Assumption \ref{ass-gn} ($p \in [0,1),q>1, r,s\geq 1$).}
	\label{tb-g}
	\centering
	\small
	\begin{tabular}{p{2cm}p{2.3cm}p{3cm}}
		\hline
		$g(y)$     &  $g_p(y)$ &  $g_q(y)$ \\
		\hline
		\multirow{2}{*}{$ g_p + g_q $}  &  $\partial \|y\|_r$ or $\frac{y}{\|y\|_r^{1-p}}$;    & $y\|y\|_r^{q-1}$;\\
				 &  $\text{sign}(y)\odot |y|^{p}$   &  $\text{sign}(y)\odot |y|^{q}$\\ 
		$\text{sign}(y)\odot e^{|y|}$  &  $\text{sign}(y)$  & $\text{sign}(y)\odot (e^{|y|}-1_n) $ \\
		$ e^{\|y\|_s}\partial \|y\|_r $  &  $\partial \|y\|_r$  & $(e^{\|y\|_s}-1) \partial \|y\|_r $ \\
		\hline
	\end{tabular}
\end{table}

\begin{thm}\label{thm1}
	Let Assumptions \ref{ass-f}, \ref{ass-f2} and \ref{ass-gn} hold. With $u(x(t))= g(\nabla f(x))$ and $d \equiv0$, any Filippov solution $x(t)$ of the nominal dynamics \eqref{pro-u} will converge to $\mathbb{X}^*$ in fixed time bounded by 
	\begin{align}\label{Tf-nominal}
		T(x_0) \leq \frac{1}{\mu\sigma (1-p) } + \frac{1}{\mu \rho(q-1)}.
	\end{align}
\end{thm}
\begin{proof}
	See Appendix B.
\end{proof}

Particularly, when considering $g(y)=\frac{ye^{\|y\|_2}}{\|y\|_2}$, we obtain the following interesting result.
\begin{prop}\label{prop-ge}
	Let Assumptions \ref{ass-f} and \ref{ass-f2} hold. For the nominal system \eqref{pro-u} with $u(x(t))= g(\nabla f(x))$ and $g(y)=\alpha \frac{ye^{\|y\|_2}}{\|y\|_2}$ (resp. $\alpha\text{sign}(y)\odot e^{|y|}$), $\alpha>0$, any solution $x(t)$ will converge to $\mathbb{X}^*$ in fixed time globally bounded by $\frac{1}{\alpha\mu}$ (resp. $\frac{n}{\alpha\mu}$). 
\end{prop}
\begin{proof}
	See Appendix C.
\end{proof}

For the completeness of this work, we provide a general condition for the FT stability of \eqref{pro-u} with respective to $\mathbb{X}^*$ for the nominal system \eqref{pro-u}. 
\begin{prop}\label{prop-FT}
	Let Assumptions \ref{ass-f} and \ref{ass-f2} hold and $g \in \mathcal{F}_p^{\sigma}$ with $p \in [0 \ 1)$ and $\sigma>0$. For the nominal system \eqref{pro-u} with $u(x)=g(\nabla f(x))$, any solution $x(t)$ starting from $x_0$ will converge to $\mathbb{X}^*$ in finite time 
	\begin{align}\label{Tf2-gn}
		T(x_0) \leq \frac{(2\mu(f(x_0)-f^*))^{\frac{1-p}{2}}}{\mu \sigma (1-p)}.
	\end{align}
\end{prop}
\begin{proof}
	With the same Lyapunov function as that of Theorem \ref{thm1}, we obtain that 
	\begin{align*}
		\frac{d}{dt}V(x(t)) = -2\mu \sigma y^Tg(y) \leq  -2\mu \sigma \|y\|_2^{p+1}=  -2\mu \sigma  V^{\frac{p+1}{2}},
	\end{align*}
	which indicates that $x(t)$ will converge to $\mathbb{X}^*$ in finite time bounded by \eqref{Tf2-gn} due to Lemma \ref{lem-fixed}. 
\end{proof}
%
%
%

Then, a unified method for designing RFxTGFs with disturbance rejections will be given for the perturbed system \eqref{pro-u} under Assumption \ref{ass-d}. Practically, the disturbance $d(t)$ can arise in an internal or external environment. On one hand, it may represent the necessary errors when calculating $u(x(t))$ related to $\nabla f(x)$ in the internal system and the distance $\|x-[x]^*\|_2$ can be used to measure $\|\nabla f(x)\|_2$ because of the inequality \eqref{QG2}. On the other hand, the disturbance can emerge in a noisy environment, such as when solving a data-driven learning task where only a noisy estimate of the gradient is available \cite{Budhraja2021}. Besides, it is more practical to consider both the bounded and diminishing disturbances in Assumption \ref{ass-d} compared with \cite{Budhraja2021} where only the latter case can be rejected by an FxTGF with continuous right-hand side. One can also replace the term $\|x-[x]^*\|_2$ with $\|\nabla f(x)\|_2$. To address non-diminishing but bounded disturbances, the discontinuous dynamics is used here coming from sliding mode control, see \cite{Shi2020TAC,Cortes2008,Yu2021}.

\begin{ass}\label{ass-d}
	The disturbance vector $d(x,t)$ is measurable and satisfies that there exist $\epsilon > 0$ and $\bar{d}> 0$ such that $\|d(x,t)\|_2\leq \epsilon \|x-[x]^*\|_2 +\bar{d}$ for any $x\in \mathbb{R}^n$.
\end{ass}

For the purpose of disturbance rejection, we suppose that Assumption \ref{ass-gn} holds with $p=0$, i.e., 
\begin{align}\label{prot-gd}
	g(\nabla f(x)) = g_0(\nabla f(x))+ g_q(\nabla f(x)),
\end{align}
where $g_0$ is a nonsmooth function (e.g., $\text{sign}(y)$, $\frac{y}{\|y\|_2}$) and $g_q \in \mathcal{F}_q^{\rho}$ with $q>1$ is used to drive $x$ into a bounded ball within a fixed time. Moreover, $g_0$ is used for disturbance rejection concerning the bounded part and together with $g_q$ for compensating the vanishing term as reflected by the following result.

\begin{thm}\label{thm-gnd}
	Let Assumptions \ref{ass-f}-\ref{ass-gn} hold with $p=0$ and let $u(x(t))= g(\nabla f(x))$. If $\sigma > \bar{d}+ \frac{\epsilon}{2\sqrt{\mu}}$ and $\rho\geq \frac{\epsilon}{2\sqrt{\mu}}$, any Filippov solution $x(t)$ of the dynamics \eqref{pro-u} will converge to $\mathbb{X}^*$ in fixed time bounded by 
	\begin{align}\label{Tf}
		T(x_0) \leq  \frac{1}{\mu k_1} + \frac{1}{\mu k_2(q-1)},
	\end{align}
	where $k_1 =\sigma -\bar{d}- \frac{\epsilon}{2\sqrt{\mu}} $ and $k_2=\rho-\frac{\epsilon}{2\sqrt{\mu}}$. 
\end{thm}
\begin{proof}
	See Appendix D.
\end{proof}

\begin{remark}
	Note that the first-order FxTGF provided in \cite{Garg2021tac} is a special case that satisfies Assumption \ref{ass-gn}, i.e., 
	\begin{align}\label{g-Garg}
		g(y) =c_1\frac{y}{\|y\|^{\frac{p_1-2}{p_1-1}}}+ c_2\frac{y}{\|y\|^{\frac{p_2-2}{p_2-1}}}
	\end{align}
	with $c_1,c_2 >0, p_1>2$ and $1<p_2<2$. Moreover, \cite{Budhraja2021} further considers the robustness of \eqref{g-Garg} in presence of the state-related vanishing disturbance. However, \eqref{g-Garg} can not be used for rejecting non-vanishing disturbances. Compared with these works, Assumption \ref{ass-gn} provides a unified framework for designing FxTGFs including both continuous and discontinuous flows, and the provided discontinuous FxTGFs can be used to reject both bounded and vanishing disturbances, as shown in Theorem \ref{thm-gnd}. Significantly, it can generate many novel protocols such that $g(y)=\alpha \frac{ye^{\|y\|_2}}{\|y\|_2}, \alpha>0$ analyzed in Proposition \ref{prop-ge}, for which the FxT bound is unrelated with the parameters $p,q$ as in \eqref{Tf-nominal} and can be prespecified by setting $\alpha = \frac{1}{T\mu}$. Moreover, Assumption \ref{ass-gn} can be further relaxed as there exist $p \in [0,1)$ and $q>1$ such that 
	\begin{align*}
		g^T(y) y \geq   \left\{\begin{aligned}
			\sigma \|y\|_2^{1+p}, &\ \|y\|\leq c, \\
			\rho\|y\|_2^{1+q}, & \ \|y\|> c.
		\end{aligned}\right.
	\end{align*}
	for some scalar $c>0$. Similarly, one can derive an FxT bound (related to $c$) regardless of initial states. 
	
\end{remark}


\subsection{Newton-Based FxTGFs and Its Online Implementation}
\begin{ass}\label{ass-f3}
	The cost function $f(x)$ is strictly convex and twice continuously differentiable with $\nabla^2 f(x)\succ 0$. 
\end{ass}

In this subsection, we will consider a Newton-based FxTGF (NFxTGF) under Assumption \ref{ass-f3} as follows
\begin{align}\label{prot-gn2}
\begin{split}
	u(x(t)) &= (\nabla^2 f(x))^{-1} g(\nabla f(x)),   
\end{split}
\end{align}
where $g$ satisfies Assumption \ref{ass-gn}. 
Protocol \eqref{prot-gn2} is also an extension of the FxT Newton’s method provided in \cite{Garg2021tac}. 
Differently from the first-order FxTGFs, $f(x)$ does not need to satisfy PL inequality \eqref{PL} but is required to be strictly convex to make \eqref{prot-gn2} well-defined, e.g., $f(x)=x+e^{-x}, x\in \mathbb{R}$. Note that the fact $f(x)$ is strictly convex and twice continuously differentiable may not guarantee that $\nabla^2 f(x)\succ 0$, e.g., $f(x)=x^4$. Thus, $\nabla^2 f(x)\succ 0$ is required in Assumption \ref{ass-f3}. In this case, there exists a unique optimal solution $x^*$, i.e., $\mathbb{X}^*=\{x^*\}$. Then, we show that the nominal dynamics \eqref{pro-u} with \eqref{prot-gn2} is also a class of FxTGFs.
\begin{thm}\label{thm-gn2}
	Let Assumptions \ref{ass-f}, \ref{ass-gn} and \ref{ass-f3} hold. With the protocol \eqref{prot-gn2} and $d\equiv0$, any Filippov solution $x(t)$ of the nominal dynamics \eqref{pro-u} starting from any initial state $x_0 \in \mathbb{R}^n$ will converge to $x^*$ in fixed time bounded by 
	\begin{align}\label{Tf-nominal2}
		T(x_0) \leq \frac{1}{\sigma (1-p) } + \frac{1}{\rho(q-1)}.
	\end{align}
Moreover, $u(x(t))$ is m.l.b. over $[0 \ T(x_0)]$. 
\end{thm}
\begin{proof}
	See Appendix E.
\end{proof}

Invoking the proofs of Propositions \ref{prop-ge} and considering the Lyapunov function $V(x(t)) = \|\nabla f(x)\|_2^2$ as in Theorem \ref{thm-gn2}, one can deduce the following results. 
\begin{prop}
	Let Assumptions \ref{ass-f} and \ref{ass-f3} hold. For the nominal system \eqref{pro-u} with \eqref{prot-gn2} and $g(y)=\alpha \frac{ye^{\|y\|_2}}{\|y\|_2}, \alpha>0$, any solution $x(t)$ will converge to $x^*$ in fixed time globally bounded by $\frac{1}{\alpha}$. 
\end{prop}

Actually, the protocol \eqref{prot-gn2} can be extended to solve the time-varying optimization, i.e., 
\begin{align}\label{CO-central}
	\min_{x\in \mathbb{R}^n} f(x,t).
\end{align}
where $f(x,t)$ is supposed to satisfy Assumptions \ref{ass-f} and \ref{ass-f3} w.r.t. $x$ for each $t$ and is differential w.r.t. $t$. Let $x^*(t)$ be the continuous optimal trajectory that solves \eqref{CO-central}. The modified FxTGFs for solving \eqref{CO-central} is given as follows
\begin{align}\label{prot-TV}
	\dot{x} =-\nabla_{xx}^{-1}f(x,t)(g(\nabla f(x,t))+\frac{\partial \nabla f(x,t)}{\partial t}).
\end{align}
\begin{prop}\label{prop-cent}
	Suppose that $f(x,t)$ satisfies Assumptions \ref{ass-f} and \ref{ass-f3} w.r.t. $x$ for each $t$ and Assumption \ref{ass-gn} holds. Then, with \eqref{prot-TV}, $x(t)$ will track $x^*(t)$ in fixed time bounded by \eqref{Tf-nominal2}. 
\end{prop}
\begin{proof}
	Considering the dynamics \eqref{prot-TV}, one can derive that 
	\begin{align*}
		\frac{d}{dt} \nabla f(x,t) = \nabla_{xx}f(x,t)\dot{x} + \frac{\partial \nabla f(x,t)}{\partial t} = -g( \nabla f(x,t)).
	\end{align*}
    Then, based on Lemma \ref{lem-fixed}, one can show that $\nabla f(x,t)$ converges to zero in fixed time bouned by \eqref{Tf-nominal2}.  
\end{proof}

\begin{remark}
Compared with the first-order methods involving only gradient information, Assumption \ref{ass-f3} is required for NFxTGFs to make the dynamics well-defined. Moreover, $f$ does not need to satisfy PL inequality for \eqref{prot-gn2}. The convergence time for the first-order FxTGFs relies on the parameter $\mu$ to the PL inequality while that of the second-order one is globally bounded by a common value regardless of the cost functions. However, the inverse of Hessian is required for the second-order method, resulting in high computation complexity. In addition, the second-order methods may be sensitive to the numerical errors while some well-designed first-order FxTGFs are more robust for disturbance rejections when external disturbance is imposed.
\end{remark}


\subsection{PFxTGFs for Constrained Optimization with Equations}\label{sub-COLE}
In this subsection, based on the first-order FxTGFs, a class of PFxTGFs will be provided for solving the following linear equality constrained optimization:
\begin{align}\label{app-linear}
	&\min_{x\in \mathbb{R}^{n}} f(x), \ \text{s.t.} \ Ax=b, 
\end{align}
where $A\in \mathbb{R}^{m\times n}$ and $b\in \mathbb{R}^m$. 
Suppose that the optimal solution set $\mathbb{X}^*\neq \emptyset$ and $x$ evolves in a time-invariant feasible set given by $\mathbb{X}_{f} = x_b + \mathcal{N}(A) = \{x=x_b+P\hat{x}: \hat{x} \in \mathbb{R}^{n} \}$, where $x_b$ satisfies $Ax_b=b$ and $P$ is a projection matrix onto $\mathcal{N}(A)$, i.e., $\mathcal{R}(P)=\mathcal{N}(A)$. Here, we do not require that $P$ is idempotent. Note that the matrix $P$ is not unique. Introduce an auxiliary function $\hat{f}(\hat{x}) = f(x_b+P \hat{x})$. To use the previous theoretical results, we should make sure that $\hat{f}(\hat{x})$ satisfies Assumption \ref{ass-f2}, i.e., PL inequality. According to \cite{Karimi2016}, when $f(x)$ is $\mu-$strongly convex, for any matrix $P\in \mathbb{R}^{n\times n}$, $\hat{f}(\hat{x})$ satisfies PL inequality:
\begin{align}\label{app-PL0}
	\frac{1}{2}\|\nabla \hat{f}(\hat{x})\|_2^2 \geq \mu \lambda_{2}(P^TP)(\hat{f}(\hat{x})-\hat{f}^*).
\end{align} 
Then, for $\hat{f}(\hat{x})$, a general FxT convergent algorithm can be designed as
\begin{align}
	\dot{\hat{x}} = -g(\nabla \hat{f}(\hat{x})) = -g(P^T \nabla f(x_b+P\hat{x}))
\end{align}
when $g$ satisfies Assumption \ref{ass-gn} and $f$ is $\mu-$strongly convex. With $x =x_b+P\hat{x}$, one can obtain the dynamics of $x$:
\begin{align}\label{prot-linear}
	\dot{x} = -P\cdot g(P^T \nabla f(x)).
\end{align}
Then, it is sufficient to make $x_0 \in \mathbb{X}_f$, i.e., $Ax_0=b$, and thus $\mathbb{X}_f$ will be time-invariant since $A\dot{x} \equiv 0$. To summarize, we give the following result.
\begin{prop}\label{lem-proj}
	Let Assumption \ref{ass-gn} hold and $f$ be $\mu-$strongly convex. If $\mathcal{R}(P)=\mathcal{N}(A)$ and $x_0 \in \mathbb{X}_f$, the dynamics of \eqref{prot-linear} will evolve in $\mathbb{X}_f$ and converge to $\mathbb{X}^*$ in fixed time bounded by
	\begin{align}\label{Tf-linear}
		T(x_0) \leq \frac{1}{\mu \lambda_{2}(P^TP) }(\frac{1}{\sigma (1-p) } + \frac{1}{\rho (q-1)}).
	\end{align}
\end{prop}

To obtain an FxT protocol with free initialization, one can first get a solution to $Ax=b$ in fixed time. For this purpose, the following algorithm is considered
\begin{align}\label{cent-fixed2}
	\dot{x} = -A^T\hat{g}(Ax-b).
\end{align}
We show the FxT convergence of \eqref{cent-fixed2} with $\hat{g}(y) \in \text{span}\{y\}$ in Proposition \ref{lem-cent}. Moreover, a disturbance rejection version of \eqref{cent-fixed2} can also be obtained by further referring to the conditions of Theorem \ref{thm-gnd}.  
\begin{prop}\label{lem-cent}
	When the function $\hat{g}$ satisfies Assumption \ref{ass-gn} and $\hat{g}(y) \in \text{span}\{y\}$, the dynamics of \eqref{cent-fixed2} will converge to the solution of $Ax=b$ in fixed time bounded by 
	\begin{align}\label{Tf-y}
		T(x_0) \leq \frac{1}{\sigma\lambda_2(AA^T) (1-p) } + \frac{1}{\rho\lambda_2(AA^T)(q-1)}.
	\end{align}
\end{prop}
\begin{proof}
	Denote $y= Ax-b$. Since the solution to $Ax=b$ exists, then $b\in \mathcal{R}(A)$ and $y\in  \mathcal{R}(A)$. With \eqref{cent-fixed2}, we have 
	\begin{align}\label{deri-y}
		\dot{y} = -AA^T\hat{g}(y). 
	\end{align}
	By Theorem \ref{thm1}, it is sufficient to show that $\tilde{g} = AA^T\hat{g}$ satisfies Assumption \ref{ass-gn} for any $y \in \mathcal{R}(A)$, which holds since $y\tilde{g}(y) \geq \lambda_2(AA^T) y^T\hat{g}(y)$ by the assumptions on $\hat{g}$. Then, the FxT bound can be deduced for $f(y) = \frac{1}{2}y^Ty$ according to \eqref{Tf-nominal}.
\end{proof}

Then, combining \eqref{prot-linear} with \eqref{cent-fixed2}, a class of PFxTGFs with free initialization can be designed by
\begin{align}\label{prot-free}
	\dot{x} = -P\cdot g(P^T \nabla f(x))-A^T\hat{g}(Ax-b).
\end{align}
The FxT convergence is shown in the following result. 
\begin{thm}
	Let $g,\hat{g}$ satisfy Assumption \ref{ass-gn}, $\hat{g}(y) \in \text{span}\{y\}$, and $f$ be $\mu-$strongly convex. If $\mathcal{R}(P)=\mathcal{N}(A)$, the dynamics of \eqref{prot-free} will converge to $\mathbb{X}^*$ in fixed time globally.
\end{thm}
\begin{proof}
	The result can be shown in two stages. First, let $y=Ax-b$. Since $AP=0$, \eqref{prot-free} reduces to \eqref{deri-y}. From Proposition \ref{lem-cent}, $x(t)$ will converge to $\mathbb{X}_f$ in fixed time and stay therein. By the assumptions on $\hat{g}$, we have $\hat{g}(0)=0$. Thereafter, the dynamics \eqref{prot-free} will be dominated by \eqref{prot-linear}. Then, by Proposition \ref{lem-proj}, $x(t)$ with converge to $\mathbb{X}^*$ in fixed time.
\end{proof}

In order to decrease the convergence time $T(x_0)$ in \eqref{Tf-linear}, it is sufficient to make $\lambda_{2}(P^TP)$ as large as possible. Without loss of generality, for a fair discussion, suppose that $P \in \mathbb{P}=\{P\in \mathbb{R}^{n\times n}: \|P\|_2\leq 1, \mathcal{R}(P)=\mathcal{N}(A)\}$. As $\lambda_{2}(P^TP)\leq \lambda_{\max}(P^TP)= \|P\|_2^2 \leq 1$, $T(x_0)$ will be minimized if $\lambda_{2}(P^TP)= \|P\|_2=1$, i.e., all the nonzero singular values of $P$ are the same. One choice for such matrix can be the orthogonal projection matrix onto $\mathcal{N}(A)$, denoted by $P_A=I-\widetilde{A}(\widetilde{A}^T\widetilde{A})^{-1}\widetilde{A}^T$, where $\widetilde{A}\in \mathbb{R}^{r\times n}$ is the submatrix of $A$ and satisfies $\text{rank}(A)=\text{rank}(\widetilde{A})=r$. Moreover, $P_A$ satisfies $P_A=P_A^T$ and $P_A^2 = P_A$. When the function $g$ in \eqref{prot-linear} satisfies $g(y) \in \text{span}\{y\}$,  \eqref{prot-linear} can be further simplified as $\dot{x} = -g(P_A \nabla f(x))$ due to $P_A^2=P_A$. 

\begin{remark}
Specifically, a class of FT solvers can be derived when $g \in \mathcal{F}_p^{\sigma}$ with $p \in [0 \ 1)$. For example, a special nonsmooth case with $g(y)=\text{sign}(y)$ is considered in \cite{Chen2020sign} for solving \eqref{app-linear} in finite time with the assumptions that $f(x)$ is convex, twice continuously differential and $A$ is of full row rank. In \cite{Garg2021tac}, for the equation-constrained optimization with the assumption that $A$ has full row rank, it is first transformed into a dual problem and then solved by the dual and primal FxTGFs, where two FxTGFs are used sequentially and the conjugate function of $f$ is required to be calculated. The same transformation is used in \cite{Guo2022CAA}. In our work, only one PFxTGF evolving the primal state is used and the matrix $A$ does not need to be of full row rank. Besides, the provided unified approach includes both discontinuous and continuous FT/FxT protocols.
\end{remark}

\subsection{FxTPGFs for Nonsmooth Composite Optimization}
Consider the following nonsmooth composite optimization problem (COP)
\begin{align}\label{COP}
	\min_{x\in \mathbb{R}^n} f(x)+h(x)
\end{align}
where $f:\mathbb{R}^n\rightarrow \mathbb{R}$ is twice continuously differentiable with Lipschitz continuous gradient, i.e., 
\begin{align}\label{f-Lip}
	\|\nabla f(x)-\nabla f(z)\|_2 \leq L_f \|x-z\|_2, \forall x,z
\end{align}
for some $L_f>0$, and $h: \mathbb{R}^n\rightarrow \mathbb{R} \cup \{+\infty\}$ is a proper, lower semi-continuous convex function. Note that if $f$ satisfies both \eqref{QG2} and \eqref{f-Lip}, one can conclude that $L_f \geq 2\mu$ by choosing $z=[x]^*$.

Denote the set of optimal solutions to \eqref{COP} as $\mathbb{X}^*$. COP \eqref{COP} can be further transformed into a constrained optimization when $h(x) =\delta_{\mathcal{X}}(x)$, where $\delta_{\mathcal{X}}(x)$ is an indicator function for a given nonempty, closed convex set $\mathcal{X}$, defined by 
\begin{equation*}
\delta_{\mathcal{X}}(x) =\left\{\begin{aligned}
		0, & \quad x\in \mathcal{X}, \\
		+\infty, & \quad x\notin \mathcal{X}.
	\end{aligned}\right.
\end{equation*}
The proximal operator of $h$ associated with a parameter $\lambda>0$ is defined as
\begin{align}\label{prox}
	\text{prox}_{\lambda h}(x) = \arg\min_{z \in \mathbb{R}^n} (h(z) + \frac{1}{2\lambda}\|z-x\|_2^2). 
\end{align}
The Moreau envelope associated with \eqref{prox} denotes its optimal value, i.e.,
\begin{align}\label{Moreau}
M_{\lambda h}(x) = h(\text{prox}_{\lambda h}(x)) + \frac{1}{2\lambda}\|\text{prox}_{\lambda h}(x)-x\|_2^2. 
\end{align}
According to \cite{Parikh2013}, $M_{\lambda h}(x)$ is essentially a smoothed
form of $h$ even if $h$ is not, and is continuously differentiable with 
$\nabla M_{\lambda h}(x) = \frac{1}{\lambda}(x-\text{prox}_{\lambda h}(x))$. The forward–backward (FB) envelope is
given by \cite{Themelis2018}
\begin{align}
	\mathcal{F}_{\lambda}(x) = f(x) + M_{\lambda h}(x-\lambda \nabla f(x))- \frac{\lambda}{2}\|\nabla f(x)\|_2^2,
\end{align}
which is continuously differentiable with 
\begin{align}
	\nabla \mathcal{F}_{\lambda}(x) = (I-\lambda \nabla^2f(x))\mathcal{H}_{\lambda}(x)
\end{align}
and $\mathcal{H}_{\lambda}(x) = \frac{1}{\lambda}(x-\text{prox}_{\lambda h}(x-\lambda \nabla f(x)))$ when $f$ is twice continuously differentiable. For the nonsmooth COP \eqref{COP}, the proximal PL inequality holds with $\lambda \in (0, \frac{1}{L_f})$ if there exists $\mu>0$ such that 
\begin{align}\label{PPL}
\frac{1}{2}	\|\mathcal{H}_{\lambda}(x)\|_2^2 \geq \mu(\mathcal{F}_{\lambda}(x)-\mathcal{F}_{\lambda}^*)
\end{align}
where $\mathcal{F}_{\lambda}^*=\arg\min_{x}\mathcal{F}_{\lambda}(x)$ which shares the same minimizers and optimal vaules with COP \eqref{COP} \cite{Hassan2021AUTO}. According to \cite{Karimi2016,Hassan2021AUTO}, \eqref{PPL} is satisfied when either of the following cases holds:
\begin{enumerate}
	\item $f$ satisfies PL inequality and $h$ is constant;
	\item $f$ is strongly convex;
	\item $f$ has the form $f(x)=\phi(Ax)$ for a strongly convex function $\phi$ and matrix $A$, while $g(x)=\delta_{\mathcal{X}}(x)$ for a convex set $\mathcal{X}$;
	\item $f$ is convex and $f+g$ satisfies the quadratic growth property like \eqref{QG}, e.g., $\|Ax-b\|_2^2/2 + \gamma \|x\|_1$ for any matrix $A$ and $\gamma>0$.
\end{enumerate}
The case 1) has been investigated in the previous parts. For case 3), $f$ (also $f+h$) is not necessary to be strongly convex. Under the proximal PL inequality condition \eqref{PPL}, for solving COP \eqref{COP}, an FxTPGF can be obtained as 
\begin{align}\label{u-cop1}
	\dot{x} = -\kappa_p \frac{\mathcal{H}_{\lambda}(x)}{\|\mathcal{H}_{\lambda}(x)\|_2^{1-p}} -\kappa_q \|\mathcal{H}_{\lambda}(x)\|_2^{q-1} \mathcal{H}_{\lambda}(x),
\end{align}
where $p \in [0, 1), q>1$ and $\kappa_p, \kappa_q>0$. 
\begin{thm}\label{thm-COP}
	For the nonsmooth COP \eqref{COP}, suppose that the proximal PL inequality \eqref{PPL} holds with $\lambda \in (0, \frac{1}{L_f})$. With the FxTPGF \eqref{u-cop1}, the trajectory $x(t)$ will converge to $\mathbb{X}^*$ in fixed time bounded by 
	\begin{align}\label{COP-Tf0}
		T(x_0) \leq \frac{1}{\mu(1-\lambda L_f)}(\frac{1}{\kappa_p (1-p) } + \frac{1}{\kappa_q (q-1)}).
	\end{align}
\end{thm} 
\begin{proof}
	Consider the candidate Lyapunov function $\widehat{V}(x(t)) = 2\mu(\mathcal{F}_{\lambda}(x(t))-\mathcal{F}_{\lambda}^*)$. By Lemma \ref{lem-non}, the time derivative of $\widehat{V}(x(t))$ can be calculated as
	\begin{align*}
		\frac{d}{dt}\widehat{V}(x) &= -\langle 2\mu \nabla \mathcal{F}_{\lambda}(x),\dot{x}\rangle \\
		&\leq -\langle2\mu (I-\lambda \nabla^2f(x))\mathcal{H}_{\lambda}(x), \dot{x}\rangle  \\
		& \leq -2\mu(1-\lambda L_f) (\kappa_p\|\mathcal{H}_{\lambda}(x)\|_2^{p+1} +\kappa_q \|\mathcal{H}_{\lambda}(x)\|_2^{q+1})	\\
		& \leq -2\mu (1-\lambda L_f)(\kappa_p (\widehat{V}(x))^{\frac{p+1}{2}}+\kappa_q (\widehat{V}(x))^{\frac{q+1}{2}})
	\end{align*}
where the last inequality is due to \eqref{PPL}. Then, by Lemma \ref{lem-fixed}, $\widehat{V}(x(t))$ will converge to zero in fixed time bounded by \eqref{COP-Tf0}, i.e., the solution $x(t)$ of \eqref{u-cop1} will converge to $\mathbb{X}^*$ in the same fixed time. 
\end{proof}
\begin{remark}\label{rmk4}
The same algorithm of the form \eqref{u-cop1} but with a different analysis method is considered in \cite{Garg2019MVI}, where the function $f$ is required to be $\mu$-strongly convex and satisfy \eqref{f-Lip}.  
Compared with \cite{Garg2019MVI}, the condition for the parameter $\gamma$ is relaxed from $(0, \frac{2\mu}{L_f^2})$ to $(0, \frac{1}{L_f})$ as $2\mu\leq L_f$. Besides, in \cite{Garg2019MVI}, $p, q$ are constrained respectively in $(1-\varepsilon(c), 1) \cap (0, 1)$ and $(1, 1+\varepsilon(c))$, where $\varepsilon(c)$ is related to parameters $\mu, \gamma$ and $L_f$, limiting its applications. However, in our results, $f$ is not required to be strongly convex and $p, q$ can be chosen freely in $(0,1)$ and $[1 \ \infty)$, respectively, based on a considerably simplified proof from a novel perspective with the proximal PL inequality. Moreover, in \cite{Ju2021Tcyber}, a novel FxT convergent forward-backward-forward neurodynamics is provided for solving mixed variational inequalities under the assumption that $\nabla f$ is $h$-strongly pesudomonotone, i.e., there exists $\mu>0$ such that 
\begin{align}
	&\langle \nabla f(x),y-x \rangle +h(y)-h(x)\geq 0 \Rightarrow \label{pm-1}\\
	 & \langle \nabla f(y),y-x \rangle +h(y)-h(x)\geq \mu \|y-z\|_2^2 \label{pm-2}
\end{align}
for any $x, y \in \mathbb{R}^n$. However, when $\eqref{PPL}$ holds, $\nabla f$ may not be $h$-strongly pesudomonotone. For example, let $f(x)=\|Ax\|_2^2$ and $h(x)=\delta_{\mathcal{X}}(x)$ with $A \in \mathbb{R}^{n\times n}$ and $\mathcal{X}$ being a convex set. When $\text{rank}(A)<n$, for any $x,y \in \mathcal{N}(A)\cap \mathcal{X}$ and $x\neq y$, \eqref{pm-1} is satisfied directly. However, there exists no constant $\mu>0$ such that \eqref{pm-2} holds. 
\end{remark}

\subsection{Regret Analysis}
In this subsection, to evaluate the effectiveness of different algorithms, we give the regret analysis based on the static regret function motivated by \cite{Budhraja2021}. For the offline optimization, the static and dynamic regret functions are the same, which can be expressed by 
\begin{align}\label{regret-fun}
	\mathcal{R}(T,x_0)  = \int_{0}^T(f(x(t))-f^*)dt, \ \forall T>0, x_0 \in \mathbb{R}^n
\end{align}
which is the accumulated difference between the cost value at $x(t)$ and the optimal one $f^*$ over the time period $[0 \ T]$. In the next, for the convenience of analysis, we give the upper bound estimation of the regret function for the nominal system \eqref{pro-u} with respect to different types of first-order protocols. By default, we denote $T(x_0)=+\infty$ for the traditional gradient descent algorithm with $g=g_1$.  
\begin{thm}\label{thm-regret}
Let Assumptions \ref{ass-f} and \ref{ass-f2} hold and consider four types of functions (i.e., $g=g_1, g_p, g_{p,q}$ and $g_e$) defined in Table \ref{tb-rb}. For the nominal system \eqref{pro-u} with $u(x)=g(\nabla f(x))$ and the initial state $x_0$, the corresponding regret bounds are shown in the second column of Table \ref{tb-rb}.
\end{thm}
\begin{proof}
	See Appendix F.
\end{proof}

\begin{remark}
In \cite{Budhraja2021}, the static regret bound is only provided for the special FxTGFs with \eqref{g-Garg}, see \eqref{Rfxt} in Appendix F. However, the estimation of the given bound is not discussed in detail. Here, we give the regret analysis for a wide range of FxTGFs, as shown in Table \ref{tb-rb}. Note that the regret bounds of $g_1,g_p$ and $g_e$ in Table \ref{tb-rb} are tight for some special functions, e.g., $f(x)=\frac{1}{2}\|x\|_2^2, x\in \mathbb{R}$. Moreover, from Table \ref{tb-rb}, it can be seen that $g_e$ has a fixed regret bound only related to $\mu$. For fixed $p \in [0,1)$ and $q>3$, the regret bound for $g_{p,q}$ is also constant unrelated with the initial states. However, when $q \in (1,3]$, the regret bound for $g_{p,q}$ depends on $V_0$ and could be unbounded. Particularly, for $f(x)=\frac{1}{2}x^2$, by replacing $g_{p,q}$ with $g(y) = \frac{y}{\|y\|_2^{1-p}} $ when $y \leq \sqrt{2}$, and $g(y) = \|y\|_2^{q-1} y $ when $y > \sqrt{2}$, the derived bound for $g_{p,q}$ in Table \ref{tb-rb} is also tight for the case $q=3$. Hence, we conclude that not all FxTGFs have fixed regret bounds regardless of the initial condition. 
\end{remark}
	
\begin{table}[t]
	\caption{Regret bounds for typical functions with $p \in [0,1), q>1, \alpha =\frac{p+1}{2} , \beta = \frac{q+1}{2}, a = (2\mu)^{\alpha}$, $b =(2\mu)^{\beta}$ and $V_0 = f(x_0)-f^*$.}
	\label{tb-rb}
	\centering
	\small
	\begin{tabular}{p{3.5cm}p{4cm}}
		\hline
		$g(y)$     &  Regret bound on $\mathcal{R}(T(x_0),x_0)$  \\
		\hline
		$g_{1}=y$  &  $\frac{V_0}{2\mu}$    \\
		$g_{p}=\frac{y}{\|y\|_2^{1-p}}$  &  $\frac{V_0^{2-\alpha}}{a(2-\alpha)}$   \\
		$g_{p,q}=\frac{y}{\|y\|_2^{1-p}} +\|y\|_2^{q-1} y $ & $\frac{1}{a(2-\alpha)}+ \frac{V_0^{2-\beta}}{b(2-\beta)}, \beta \in (1,2)$; $\frac{1}{a(2-\alpha)} + \frac{\text{ln}(1+V_0)}{b}, \beta =2$;  $\frac{1}{a(2-\alpha)}+ \frac{1}{b(\beta-2)}, \beta >2 $\\
		$g_{e}= \frac{ye^{\|y\|_2}}{\|y\|_2} $  &  $\frac{1}{\mu^2}$ \\
		\hline
	\end{tabular}
\end{table}

\section{Applications in Network Problems}\label{application}

In this part,  the proposed FxTGFs are applied to two typical network problems in a unified framework compared with existing works. Particularly, with special forms of the function $g$ satisfying Assumption \ref{ass-gn}, the considered problems can be solved in a distributed way. 
\subsection{Network Consensus}
Consider the network consensus problem over an undirected graph $\mathcal{G}$ composed of $N$ agents, which can be formulated as the following unconstrained optimization problem
\begin{align}\label{DO-consensus}
	\min_{x\in \mathbb{R}^N} f(x) = \frac{1}{2}x^TLx, 
\end{align}
where $L \in \mathbb{R}^{N\times N}$ is the Laplacian matrix of the network $\mathcal{G}$. Without of loss of generality, the state of the $i$-th agent is denoted by $x_i \in \mathbb{R}$. Under Assumption \ref{ass-G}, the set of solutions to \eqref{DO-consensus} is represented by $\mathbb{X}^* = \text{span}\{\bm{1}_N\}$ and the optimal value is $f^* =0$. Moreover, $f(x)$ satisfies PL inequality as
\begin{align}
	\frac{1}{2}\|\nabla f(x)\|_2^2 = \frac{1}{2}\|Lx\|_2^2 \geq \lambda_2(L)(f(x)-f^*),
\end{align}
where $\lambda_2(L)$ is the minimum nonzero eigenvalue of $L$, i.e., the algebraic connectivity of $\mathcal{G}$. 
\begin{ass}\label{ass-G}
	The undirected network $\mathcal{G}$ is connected.
\end{ass}

Based on the theoretical results for the FxTGFs without disturbances in Subsection \ref{subsection-gn}, one can obtain a unified approach for designing FxT consensus protocols as stated in Corollary \ref{cor-consensus}, which is a direct result of Theorem \ref{thm1}. The component-wisely sign-preserving function $g$ can be chosen as $g(y)= \text{sign}(y)\odot |y|^{1-p}+\text{sign}(y)\odot |y|^{q}$ with $p \in [0,1), q>1$. Furthermore, with $g(y)= \text{sign}(y)\odot e^{|y|}$, Corollary \ref{cor-consensus2} can be derived by Proposition \ref{prop-ge}. 
\begin{corollary}\label{cor-consensus}
	Let Assumption \ref{ass-G} hold and the component-wisely sign-preserving function $g$ satisfy Assumption \ref{ass-gn}. Consider the network consensus problem \eqref{DO-consensus}. By the nominal dynamics \eqref{pro-u} with the distributed protocol $u(x)=g(Lx)$, all the states will achieve consensus in a fixed time bounded by 
	\begin{align}\label{Tf-consensus}
				T(x_0) \leq \frac{1}{\lambda_2(L)\sigma (1-p) } + \frac{1}{\lambda_2(L)\rho (q-1)}.
	\end{align}
\end{corollary}

\begin{corollary}\label{cor-consensus2}
	Let Assumption \ref{ass-G} hold for the network problem \eqref{DO-consensus}. By the distributed protocol $u(x)=g(Lx)$ with $g(y)= \alpha \text{sign}(y)\odot e^{|y|}$, all the states of the nominal dynamics \eqref{pro-u} will achieve consensus in a fixed time bounded by $\frac{N}{\alpha \lambda_2(L)}$.
\end{corollary}

When there exist disturbances arising from communication noises or computation errors, one can deduce a unified scheme for designing FxT consensus protocols with disturbance rejections as indicated in Corollary \ref{cor-consensusd}, implied by Theorem \ref{thm-gnd}. 
\begin{corollary}\label{cor-consensusd}
	Let Assumptions \ref{ass-d} and \ref{ass-G} hold for the network consensus problem \eqref{DO-consensus} and the component-wisely sign-preserving function $g$ satisfy Assumption \ref{ass-gn} with $p=0$. Consider the disturbed system \eqref{pro-u} with the distributed protocol $u(x)=g(Lx)$. If $\sigma > \bar{d}+ \frac{\epsilon}{2\sqrt{\lambda_2(L)}}$ and $\rho\geq \frac{\epsilon}{2\sqrt{\lambda_2(L)}}$, all the states will achieve consensus in a fixed time bounded by
	\begin{align}\label{Tf-consensus2}
	T(x_0) = \frac{1}{\lambda_2(L) k_1} + \frac{1}{\lambda_2(L) k_2(q-1)},
\end{align}
where $k_1 =\sigma -\bar{d}- \frac{\epsilon}{2\sqrt{\lambda_2(L)}} $ and $k_2=\rho-\frac{\epsilon}{2\sqrt{\lambda_2(L)}}$. 
\end{corollary}

\begin{remark}
The provided distributed FxTGFs for solving network consensus problems belong to a class of node-based methods, including \cite{Xiao2009, Shi2020TAC, Shi2019TNSE,Zuo2014} as special cases, and can further be extended to cases on directed graphs and bipartite consensus. From the viewpoint of PL inequality, the theoretical analysis is simplified significantly. Moreover, with Corollary \ref{cor-consensusd}, it can be also used for disturbance rejection.
\end{remark}

\subsection{Solving A System of Linear Equations}
The second application is to solve a system of linear algebraic equations (i.e., $Ax=b$) over a connected and undirected graph $\mathcal{G}$, which has been extensively studied based on discrete- or continuous-time methods, see the survey paper \cite{Wand2019arc} and references therein. It is supposed that each agent $i$ only knows a few rows of matrix $A$ and vector $b$, i.e., $A_i \in \mathbb{R}^{m_i \times n}$ and $b_i \in \mathbb{R}^{m_i}$. Denote $x_i\in \mathbb{R}^{n}$ as the local estimation of the global solution. Let $\bm{x}=[x_i]_{i\in \langle N \rangle}$, $\hat{b}=[b_i]_{i\in \langle N \rangle}$, $\widehat{A}=\text{blkdiag}(A_1,\dots,A_N)$ and $\widehat{L}=L\otimes I_n$ with $\otimes$ being the Kronecker product. Suppose that the solution set $\mathbb{X}^*$ to $Ax=b$ is nonempty. Then, the problem can be formulated in a compact form as the following unconstrained optimization:
\begin{align}\label{app-equation}
	\min_{\bm{x}\in \mathbb{R}^{Nn}} f(\bm{x})=\frac{1}{2} \bm{x}^T\widehat{L}\bm{x} + \frac{\delta}{2}\|\widehat{A}\bm{x} -\hat{b}\|_2^2
\end{align}
with $\delta >0$. The first term in \eqref{app-equation} is to achieve consensus among agents, i.e., $x_i=x_j, \forall i,j \in \langle N \rangle $. Then, the set of optimal solutions to \eqref{app-equation} is $\bm{1}_N\otimes \mathbb{X}^*$ and $f^* = 0$ as $\mathbb{X}^*\neq \emptyset$ and $f\geq 0$. By introducing any $\bm{x}^* \in \bm{1}_N\otimes \mathbb{X}^*$ and setting $\bm{\hat{x}} = \bm{x}-\bm{x}^*$, using the fact $b=\widehat{A}\bm{x}^*$, the cost function can be rewritten as 
$f = \frac{1}{2}\bm{\hat{x}}^T(\widehat{L}+\delta \widehat{A}^T\widehat{A})\bm{\hat{x}}$, which satisfies the PL inequality as 
\begin{align}
	\frac{1}{2}\|\nabla f(\bm{x})\|_2^2 \geq \lambda_2(\widehat{L}+\delta \widehat{A}^T\widehat{A}) (f(\bm{x})-f^*).
\end{align}
Similarly to Corollaries \ref{cor-consensus}-\ref{cor-consensusd} obtained for the network consensus problems, with the component-wisely sign-preserving function $g$, one can derive a unified approach for designing distributed FxT convergent protocols to solve $Ax=b$ with or without disturbances, i.e., $\dot{\bm{x}} = -g(\widehat{L}\bm{x}+\delta \widehat{A}^T(\widehat{A}\bm{x}-\hat{b}))$, which extends the linear convergent algorithm \cite{Yang2015} to FxT convergent one. 

Specifically, when considering only one agent, a centralized FxT convergent protocol can be further provided with $f(x)=\frac{1}{2}\|Ax-b\|_2^2$:
\begin{align}\label{cent-fixed}
	\dot{x} = -g(A^T(Ax-b)), 
\end{align}
where $g$ satisfies Assumption \ref{ass-gn}. 

Another form of linear equation is considered when $A$ is partitioned by columns, resulting in
\begin{align}\label{eq-col}
	\sum_{i=1}^{N}A_ix_i = \sum_{i=1}^{N}b_i
\end{align}
where $A_i\in \mathbb{R}^{m\times n_i}$ and $b_i\in \mathbb{R}^m$ are only available to the $i$-th agent. To solve \eqref{eq-col} in a distributed mode, we introduce variables $\bm{y}=[y_i]_{i\in \langle N \rangle}$ and $\bm{z}=[z_i]_{i\in \langle N \rangle}$ with $y_i,z_i \in \mathbb{R}^m$ being local auxiliary variables of agent $i$. Then, with the similar representations as in \eqref{app-equation}, \eqref{eq-col} can be transformed into the following equivalent form 
\begin{align}\label{eq-col2}
	\bm{y} + \widehat{A} \bm{x} -\hat{b}=0; \ \bm{y}+\widehat{L}\bm{z} = 0.
\end{align} 
It can be further formulated as an unconstrained quadratic optimization
\begin{align}
	\min_{x\in \mathbb{R}^{Nn}} f(\bm{x},\bm{y},\bm{z})=\frac{1}{2} \|\bm{y} + \widehat{A} \bm{x} -\hat{b}\|_2^2 + \frac{1}{2} \|\bm{y}+\widehat{L}\bm{z}\|_2^2.
\end{align} 
As $f(\bm{x},\bm{y},\bm{z})$ is quadratic and there exists solution $(\bm{x}^*,\bm{y}^*,\bm{z}^*)$ to \eqref{eq-col2}, it satisfies the PL inequality. Then, with the component-wisely sign-preserving function $g=[g_x;g_y;g_z]$ satisfying Assumption \ref{ass-gn}, one can design distributed FxT convergent algorithms to solve \eqref{eq-col} with or without disturbances as
\begin{align*}
	\dot{\bm{x}} &=-g_x(\widehat{A}^T(\bm{y}+\widehat{A} \bm{x} -\hat{b})),\\
	\dot{\bm{y}} &=-g_y(2 \bm{y}+\widehat{A} \bm{x} -\hat{b} +\widehat{L}\bm{z} ),\\
	\dot{\bm{z}} &=-g_z(\widehat{L}(\bm{y} + \widehat{L}\bm{z})).
\end{align*}
Note that in the dynamics of $\bm{z}$, two rounds of communication related to the local information of $\bm{z}$ are required.


\begin{remark}
In the survey paper \cite{Wand2019arc}, various centralized and distributed algorithms in the viewpoint of discrete- and continuous-time methods have been discussed. However, most of these existing algorithms can only achieve asymptotic or exponential convergence. In \cite{Zhou2019tac}, a nonsmooth projection-based distributed solver is provided for linear equations with minimum $l_1$ in finite time with the assumption that $A$ has full row rank. With the same assumption on $A$, in the work \cite{Shi2020auto}, several projection-based FT and FxT continuous solvers are provided under different initial conditions. These FT/FxT projection-consensus algorithms restrain the local state to the set $\mathbb{X}_i=\{x: A_ix=b_i\}$ and then make all the agents reach consensus in finite time without disturbances. In \cite{Shi2020TAC}, a discontinuous networked dynamics is provided for solving $Ax=b$, where $A$ should be square and Lyapunov
diagonally (semi)stable. Differently from these works, without the projection term and any assumptions on $A$, by combining the consensus term with the local equations as in \eqref{app-equation}, we derive a unified scheme for designing centralized and distributed FxT convergent solvers to $Ax=b$ with or without disturbances for free initialization. 
\end{remark}

\section{Numerical examples}\label{Numberical}
In the following, four case studies will be performed to solve different types of optimizations by applying the previous proposed FxTGFs with/without disturbances compared with existing methods. 
\subsection{Case 1: Logistic Regression With Disturbance}
\begin{figure}[t]
	\centering
	{\includegraphics[width=.35\textwidth]{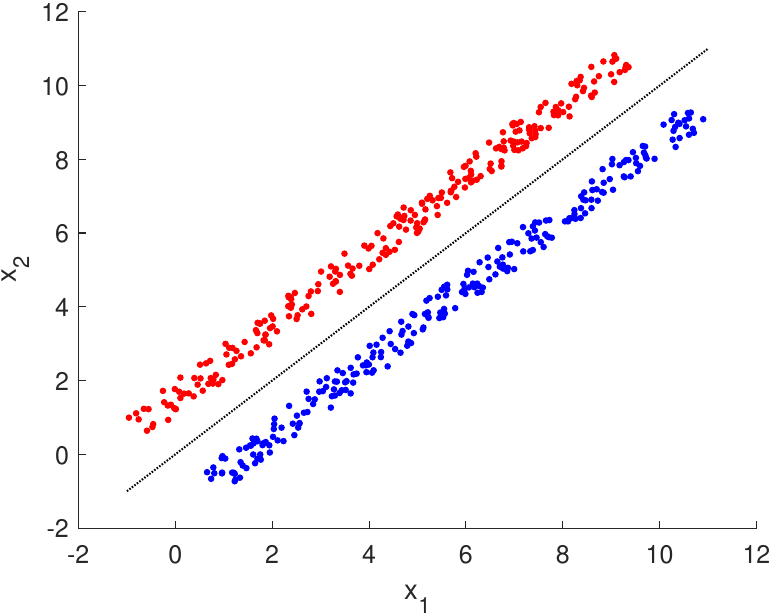}}
	\caption{Distribution of sample data $z_k$ around the line $x_1=x_2$ with red and blue dots labeled by $l_k=-1$ and $1$, respectively.}
	\label{case1-f1}
\end{figure}

\begin{figure}[t]
	\centering
	{\includegraphics[width=.35\textwidth]{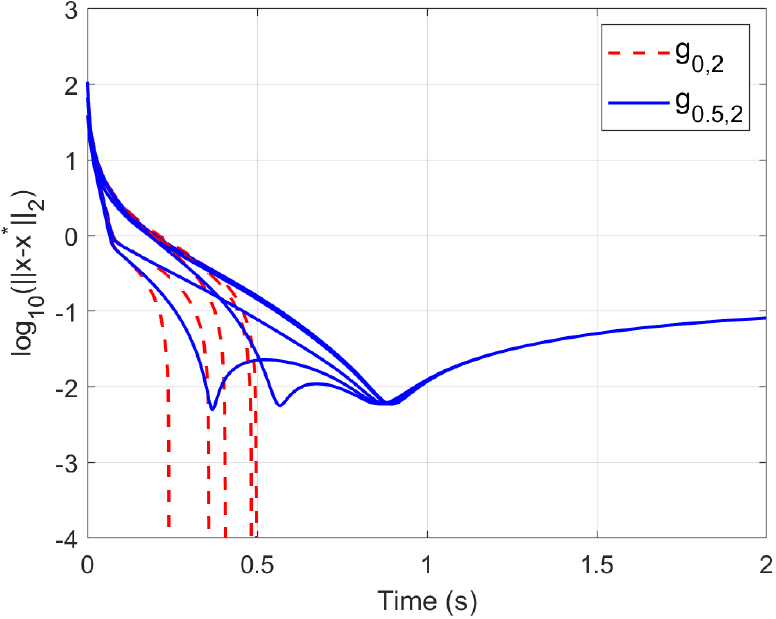}}
	\caption{Case 1: Values of $\log_{10}(\|x-x^*\|_2)$ for protocols $g_{0,2}$ (dashed lines) and $g_{0.5,2}$ (solid lines) with different initial conditions.}
	\label{case1-f2}
\end{figure}

Consider the following binary logistic regression model with $l_2$-regularization:
\begin{align}\label{LR}
	\min_{x} f(x) = \frac{1}{K}\sum_{k=1}^K\log(1+\exp(-l_kx^Tz_k)) +\frac{\beta}{2} \|x\|_2^2,
\end{align} 
where $l_k \in \{-1,1\}$ denotes the label of the $k$-th sample data $z_k$ and $\beta$ is the
regularization parameter. In the case study, sample vectors $z_k$ are randomly distributed around the line $x_1=x_2$, as shown in Fig. \ref{case1-f1}. Then, the objective of \eqref{LR} is to find a hyperplane separating the given data set $\{z_k: k\in \langle K \rangle\}$. In the simulation, $\beta=1$, $K=500$ and the disturbance in the dynamics \eqref{pro-u} is set by $d(t)=[\sin(t) \ \cos(t)]^T(1+\|x(t)-x^*\|)$ including both vanishing and bounded disturbances to satisfy Assumption \ref{ass-d}. In this case, the robust protocol $g_{0,2}$ that can reject both vanishing and bounded disturbances is compared with the protocol $g_{0.5,2}$ studied in \cite{Budhraja2021} designed only for vanishing disturbance rejection, which are represented by:
\begin{align*}
	g_{0,2}(y) =  \frac{3y}{\|y\|_2} +  3 y\|y\|_2^2; \ g_{0.5,2}(y) =  \frac{3y}{\|y\|_2^{0.5}} +  3 y\|y\|_2^2
\end{align*}
to satisfy the conditions provided in Theorem \ref{thm-gnd} and \cite[Thm. 2]{Budhraja2021}, respectively.
The initial states are randomly chosen in $[-100,100]\times [-100,100]$ and five runs are executed with different initialization. The simulation results are shown in Fig. \ref{case1-f2}. Since $f(x)$ is 1-strongly convex, it indicates that $f(x)$ satisfies PL inequality with parameter $\mu=1$. By Theorem \ref{thm-gnd}, one can calculate the theoretical time bound by $T^* =1.0667$ with $\bar{d}=1, \epsilon=1, \sigma=\rho=3$ and $q=2$. From Fig. \ref{case1-f2}, it can be seen that the state $x(t)$ with $g_{0,2}(y)$ converges to $x^*$ in fixed time bounded by $T^*$ regardless of the initial states and the disturbance. However, with the protocol $g_{0.5,2}(y)$, the state only converges to a neighborhood of $x^*$ and fluctuates therein because of non-vanishing disturbances.

\begin{figure*}[htbp]
	\centering
	\begin{minipage}[t]{0.31\textwidth}
		\centering
		\includegraphics[width=1\textwidth]{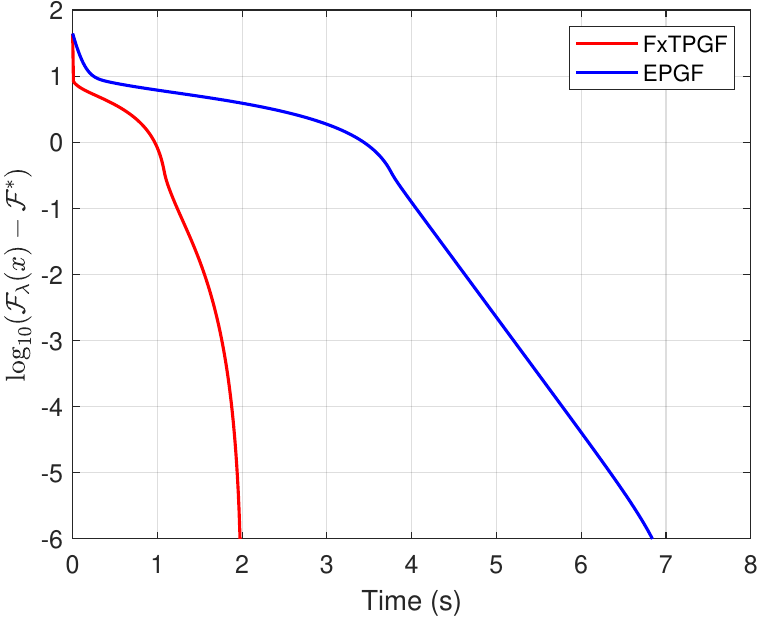}
        \caption{Case 2: Values of $E_f(t)$ for FxTPGF \eqref{u-cop1} and EPGF \cite{Hassan2021AUTO}.}
       \label{case2-FE}
	\end{minipage}
	\begin{minipage}[t]{0.31\textwidth}
		\centering
		\includegraphics[width=1\textwidth]{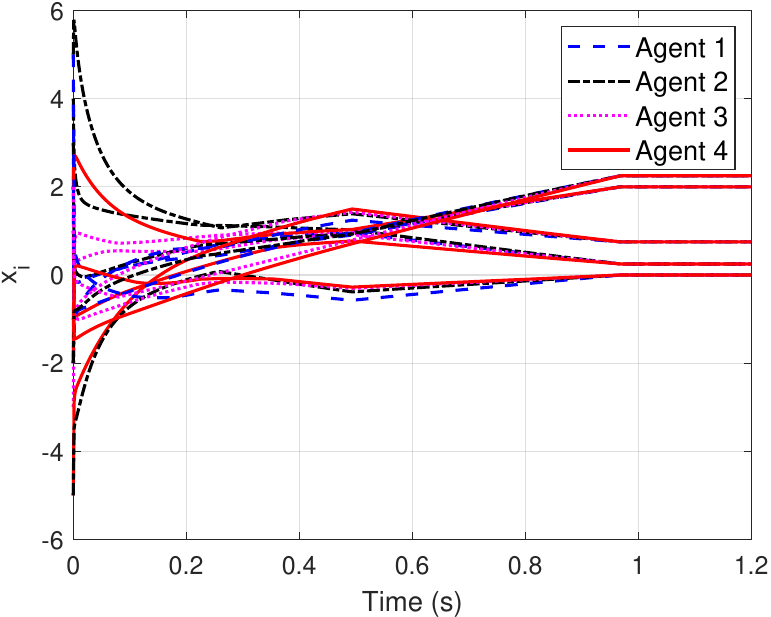}
	\caption{Case 3: Values of $\log_{10}(\|\nabla f(x)\|_2)$ for the protocols with $g_{d}$ and $g_{c1/c2}$.}
    \label{case3-f1}
	\end{minipage}
	\begin{minipage}[t]{0.31\textwidth}
		\centering
		\includegraphics[width=1\textwidth]{Fig6.pdf}
	\caption{Case 3: Agents' states by the distributed protocol with $g_{d}$.}
    \label{case3-f2}
	\end{minipage}
\end{figure*}


%

\subsection{Case 2: Constrained Lasso Problem}
Consider the following Lasso problem,
\begin{align}
	\min_{x\in \mathcal{X}} \frac{1}{2}\|Ax-b\|_2^2+\gamma \|x\|_1
\end{align}
with $\gamma=1$, $\mathcal{X} = \{x\in \mathbb{R}^4: -5\leq x_i \leq 5, i \in \langle 4\rangle \}$ and 
\begin{align}
	A &= \left[
	\begin{array}{cccc}
		1& 0& -1& 0\\
		1& 2& -1&-1 \\
		0& 0& 0&1\\
	\end{array}
	\right]; b = \left[
	\begin{array}{c}
		1 \\
		1 \\
		1\\
	\end{array}
	\right].
\end{align}
It can be transformed into COP \eqref{COP} with $f(x)=\frac{1}{2}\|Ax-b\|_2^2$ and $h(x)=\|x\|_1+ \delta_{\mathcal{X}}(x)$. According to \cite{Karimi2016}, $f$ satisfies the proximal PL inequality
\eqref{PPL}. Moreover, $f$ has Lipschitz continuous gradient with constant $L_f =7.83$. Then, we implement the provided FxTPGF \eqref{u-cop1} compared with the exponentially convergent proximal gradient flow (EPGF) (i.e., $\dot{x}=-\mathcal{H}_{\lambda}(x)$) \cite{Hassan2021AUTO} for comparison. Note that the same FxTPGF \eqref{u-cop1} is considered in \cite{Garg2019MVI} but with more strict conditions on parameters $\gamma, p, q$, and $f(x)$ needs to be strongly convex, as discussed in Remark \ref{rmk4}. Hence, in this case, only the general settings are considered for parameters and $f(x)$ is not required to be strongly convex, as stated in Theorem \ref{thm-COP}. In the simulation, we set $\lambda=0.1$, $p=0.5, q=2, \kappa_p=\kappa_q=1$ for \eqref{u-cop1}, and the initial state is randomly chosen in $\mathcal{X}$. As the optimal solution is not unique, the exponential convergence of $\|x-x^*\|_2$ may not be established even though $\mathcal{F}_{\lambda}(x)-\mathcal{F}^*$ converges to zero exponentially \cite{Hassan2021AUTO}. Here, we use the error function $E_f(t) = \log_{10}(\mathcal{F}_{\lambda}(x)-\mathcal{F}^*)$ in terms of function values to measure the convergence rate of two algorithms, where $\mathcal{F}^* = 1.25$ is the optimal objective value to the considered COP. By the simulation, the values of $E_f(t)$ are shown in Fig. \ref{case2-FE}. It can be seen that the proposed FxTPGF \eqref{u-cop1} converges to the optimal solution in a finite time less than 2s, which is much faster than EPGF.

\subsection{Case 3: Solving Linear Equations}
The third case study is to apply the proposed distributed and centralized FxT protocols for solving the linear equations $Ax=b$ with 
\begin{align*}
	A &= \left[
	\begin{array}{ccccc}
		3&4&-3&-2&-2\\
		1&-2&-4&-5&3\\
		\hdashline[2pt/2pt]              
		4&5&-2&-2&-2\\
		0&-4&4&4&4\\
		\hdashline[2pt/2pt]
		3&-4&-3&4&2\\
		\hdashline[2pt/2pt]
		5&-3&-5&-5&2\\
	\end{array}
	\right], b=\left[
	\begin{array}{c}
		2\\
		0\\
		\hdashline[2pt/2pt]
		5\\
		4\\
		\hdashline[2pt/2pt]
		-5\\
		\hdashline[2pt/2pt]
		-4\\
	\end{array}
	\right], 
\end{align*}
the integer elements of which are randomly generated from $[-5,5]$ but should satisfy $\text{rank}(A) = \text{rank}([A\ b])$ to guarantee the existence of solutions. For the distributed algorithm, the unconstrained optimization \eqref{app-equation} is considered and the protocol function is denoted by $g_d(y)$ with $y = \nabla f(\bm{x})$. The communication graph is a circle composed of $N=4$ agents and the edge weights are all ones. Two centralized algorithms given by \eqref{cent-fixed} and \eqref{cent-fixed2} are also used and the involved protocol functions $g$ are denoted by $g_{c1}(y)$ and $g_{c2}(y)$, respectively. The related functions have forms as follows 
\begin{align*}
	&g_{d}(y) =g_{c1}(y)= 3 \text{sign}(y)+  3 \text{sign}(y)\odot |y|^{1.5}; \\
	&g_{c2}(y) =  \frac{3y}{\|y\|_2} +  3 y\|y\|_2^{0.5}.
\end{align*}

Besides, the following distributed edge-based projected algorithm (EPA) proposed in \cite{Shi2020auto} is realized as a comparison 
\begin{align}
	\dot{x}_i=& - 3 P_i\sum_{j\in \mathcal{N}_i}(\text{sgn}^{0.5}(x_i-x_j)-\text{sgn}^{1.5}(x_i-x_j)) \nonumber \\ 
	&- 3 A_i^T(\text{sgn}^{0.5}(A_ix_i-b_i)+ \text{sgn}^{1.5}(A_ix_i-b_i)) \label{dprojected}
\end{align}
for any $i \in \langle 4 \rangle$, where $\text{sgn}^{\alpha}(y)=\text{sign}(y)\odot |y|^{\alpha}$, $\mathcal{N}_i$ is the neighborhood of agent $i$, and $P_i=I-A_i^T(A_iA_i^T)^{-1}A_i$ is the orthogonal projection matrix on $\mathcal{N}(A_i)$. The second term of \eqref{dprojected} is used to drive the local state to the set $\mathcal{X}_i = \{x: A_ix_i=b_i\}$. However, it cannot converge to the optimal solution in the presence of disturbances. Let the initial states of the distributed and centralized algorithms be randomly chosen in $[-5, 5]^{20}$ and $[-5, 5]^{5}$, respectively, where $\mathcal{S}^r$ denotes the $r$-ary Cartesian power of the set $S$. In addition, the disturbance function $d(t) =0.2\sin(t)$ is added in the local dynamics for all the algorithms. As the solution to the linear equations is not unique, we use the common logarithm values of the gradient norm $\|\nabla f(x)\|_2$ to measure the performance of the algorithm. Specifically, $f(\bm{x})=\frac{1}{2} \bm{x}^T\bm{L}\bm{x} + \frac{1}{2}\|\widehat{A}\bm{x} -\hat{b}\|_2^2 $ for two distributed algorithms, and $f(x) =\frac{1}{2}\|Ax-b\|_2^2 $ for the other two centralized ones. The simulation results are shown in Fig. \ref{case3-f1} and the corresponding states for the distributed algorithm are illustrated in Fig. \ref{case3-f2}. It can be seen that for the former three algorithms with discontinuous term, the states converge to the solution of linear equations in finite time in the presence of disturbances. As \eqref{dprojected} is designed for normal system absent of disturbances, the states converge to a neighborhood of the optimal one, as discussed in \cite[Lemma 2]{Shi2020auto}.

\begin{table}[t]
	\caption{Generator cost parameters and local demands.}
	\label{table1}
	\centering
	\small
	\begin{tabular}{lllll}
		\hline
		Unit     &  1 &  2& 3 & 4\\
		\hline
		$a_i\ (\$$/MW$^2$h)  &   0.001562 & 0.00194 & 0.00482 & 0.00228 \\
		$b_i\ (\$$/MWh)  &  7.92  & 7.85 & 7.97&7.48 \\
		$c_i\ (\$$/h)  &561 & 310  & 78 & 459\\
		$d_i\ $(MW)  &60 & 40  & 50 & 80\\
		\hline
	\end{tabular}
\end{table}

\subsection{Case 4: Economic Dispatch Problem}
\begin{figure}[t]
	\centering
	{\includegraphics[width=.35\textwidth]{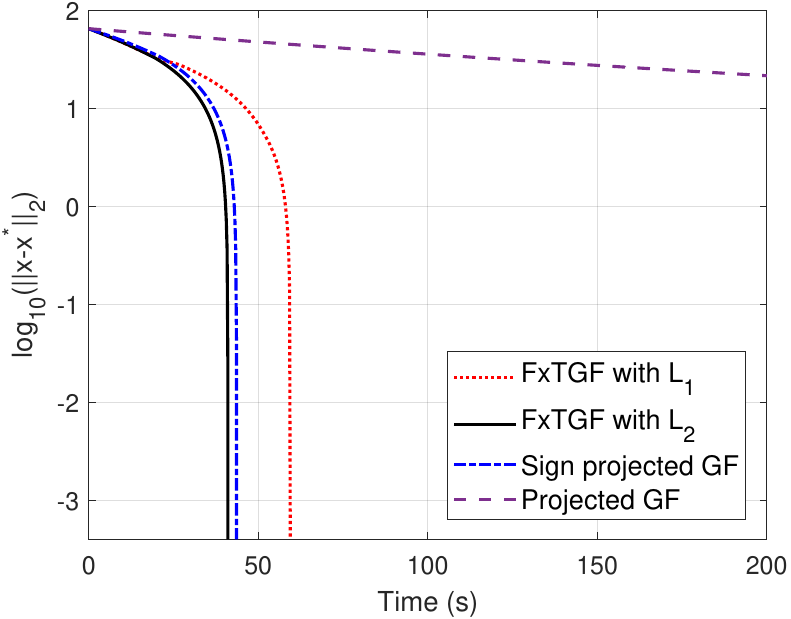}}
	\caption{Case 4: Values of $\log_{10}(\|x-x^*\|_2)$ for two FxTGFs dynamics compared with sign projected gradient flow \cite{Chen2020sign} and classic projected gradient flow \cite{Cherukuri2015}.}
	\label{case4-f1}
\end{figure}

The fourth case study is to apply \eqref{prot-linear} for solving the traditional economic dispatch in power system given as 
\begin{align}\label{ex-do2}
	\min \sum_{i=1}^4 f_i(x_i), \ s.t. \ \sum_{i=1}^4 x_i = \sum_{i=1}^4 d_i, 
\end{align}
with $d_i$ being the local demand in the $i$-th area and the local cost function $f_i(x_i) =a_i x_i^2 +b_i x_i+c_i$, whose coefficients are partially obtained from \cite{Wood2013}, as shown in Table \ref{table1}. The local information of $f_i$ and $d_i$ is only known by the $i$-th agent. 
Two projection matrices to the linear equation are chosen as 
\begin{align*}
	L_1&= \frac{1}{4}\left[
	\begin{array}{cccc}
		2 & -1 & 0 &-1\\
		-1 & 2&-1& 0 \\
		0&-1&-2&-1\\
		-1& 0 &-1& 2 \\
	\end{array}
	\right];\\
	L_2&= \frac{1}{4}\left[
	\begin{array}{cccc}
		3 & -1 & -1 &-1\\
		-1 & 3&-1& -1 \\
		-1&-1&-3 &-1\\
		-1& -1 &-1& 3 \\
	\end{array}
	\right], 
\end{align*}
which can be regarded as the Laplacian matrices of a circle undirected graph and completed graph, satisfying $\|L_1\|_2 =\|L_2\|_2 =1 $ and $\lambda_2(L_1) =0.5 < \lambda_2(L_2) =1$. In this case, two FxT dynamics \eqref{prot-linear} are executed by setting $P=L_1$ and $L_2$, respectively, and the function $g(y) = \text{sign}(y)+ \text{sign}(y)\odot |y|^{1.5}$. The proposed FxT algorithms are further compared with the sign projected GF \cite{Chen2020sign}, i.e., $\dot{x}=-L_2\text{sign}(L_2 \nabla f(x))$, and the classic projected GF--Laplacian-gradient dynamics \cite{Cherukuri2015}, i.e., $\dot{x}=-L_2 \nabla f(x)$. The latter two compared algorithms can be seen as a special case of \eqref{prot-linear} by simply setting $g = \text{sign}$ and the identity map, respectively. Note that all the algorithms are distributed since only neighbors' gradient information is involved. In the simulation, the initial states are chosen as $x_i(0)=d_i, \forall i\in \langle 4 \rangle$, which guarantees that the balance equation will not be violated along the evolution of the dynamics. The error function $E_x(t) = \log_{10}(\|x-x^*\|_2)$ is used to measure the convergence rate of all the algorithms, where $x^*=(42.87,52.56,8.71,125.86)$ is the unique optimal solution to \eqref{ex-do2}. The simulation results are shown in Fig. \ref{case4-f1}. One can see that the classic gradient projected flow has a linear convergence rate, and the other three dynamics converge to $x^*$ in finite time. Besides, the FxT dynamics with $L_2$ is faster than that with $L_1$, which verifies the analysis on the settling time at the end of subsection \ref{sub-COLE}. 

\section{Conclusion}\label{conclusion}
This work provides a unified scheme for designing fixed-time gradient flows (FxTGFs) based on nonsmooth analysis, which extends the existing finite- and fixed-time algorithms. The robustness of a class of nonsmooth FxTGFs with disturbance rejection is shown in the presence of both bounded and vanishing disturbances. PFxTGFs and FxTPGFs are provided for solving equation-constrained optimization and composite optimization, respectively, relaxing the conditions of existing results. Besides, to measure the effectiveness of the provided algorithms, the static regret analysis for several typical FxTGFs is provided in detail, and we conclude that not all FxTGFs have fixed regret bounds regardless of the initial condition. Furthermore, the proposed FxTGFs are applied to two typical network problems, for which distributed algorithms can be derived by choosing a class of component-wisely sign-preserving functions. In the future, we will focus on FxTGFs for solving constrained distributed optimization and the discretization of continuous-time FxTGFs for the practical implementation. 

\section*{Appendix}
\subsection{Proof of Lemma \ref{lem-g}}\label{proof-g}
\begin{proof}
1) Consider the function $g_0(y) = \partial \|y\|_r$ with $r\geq 1$. When $r=1$, $g_0(y)=\text{sign}(y)$ and for any $\eta \in  \mathcal{F}[\text{sign}](y)$, we have $\eta^Ty = \|y\|_1 \geq \|y\|_2$. Thus, \eqref{gp} holds with $\sigma=1$. When $q>1$, $g_0(y) = \frac{\text{sign}(y)\odot|y|^{r-1}}{\|y\|_r^{r-1}}$, which is continuous for any $y \neq 0$. And, we can derive that 
\begin{align*}
	g_0^T(y)y = \frac{\|y\|_r^r}{\|y\|_r^{r-1}} = \|y\|_r \geq \left\{\begin{aligned}
		\|y\|_2, & \quad r\in [1,2], \\
		n^{\frac{1}{r}-\frac{1}{2}} \|y\|_2, & \quad r> 2.
	\end{aligned}\right.
\end{align*} 
2) For $g_p(y) =  \frac{y}{\|y\|_r^{1-p}}, p \in (0,1)$, the second statement can be shown by  
\begin{align}
	g_p^T(y) y = \frac{\|y\|_2^2}{\|y\|_r^{1-p}} \geq \left\{\begin{aligned}
		n^{\frac{1-p}{2}-\frac{1-p}{r}}\|y\|_2^{1+p}, & \quad r\in [1,2], \\
		\|y\|_2^{1+p}, & \quad r> 2.
	\end{aligned}\right.
\end{align}

3) The assertion for $g_q(y)=y\|y\|_r^{q-1}$ is shown by the following fact
\begin{align*}
	g_q^T(y)y \geq \rho \|y\|_2^2 \|y\|_r^{q-1}\geq \left\{\begin{aligned}
		\|y\|_2^{q+1}, & \quad r\in [1,2], \\
		n^{(q-1)(\frac{1}{r}-\frac{1}{2})} \|y\|_2^{q+1}, & \quad r> 2.
	\end{aligned}\right.
\end{align*}

4) For $g(y)=\text{sign}(y)\odot |y|^{\alpha}$ with $\alpha \geq 0$, it satisfies 
\begin{align*}
	g^T(y)y =\sum_{i=1}^n |y_i|^{\alpha+1} \geq \left\{\begin{aligned}
		\|y\|_2^{\alpha+1}, & \quad \alpha \in [0,1), \\
		n^{(1-\frac{\alpha+1}{2})} \|y\|_2^{\alpha+1}, & \quad \alpha > 1.
\end{aligned}\right.
\end{align*}

5) The first part of assertion 5) holds straightforwardly due to \eqref{ge-property}. Then, it remains to show 
\begin{align*}
	\sum_{i=1}^{n} |y_i|e^{|y_i|} \geq \frac{\|y\|_2}{\sqrt{n}}e^{\frac{\|y\|_2}{\sqrt{n}}}.
\end{align*} 
It is sufficient to show that $\sum_{i=1}^{n} \frac{|y_i|^{k+1}}{k!} \geq  \frac{1}{k!}(\frac{\|y\|_2}{\sqrt{n}})^{k+1}$ for any integer $k\geq 0$, which is reduced to show $\|y\|_{k+1} \geq \frac{\|y\|_2}{\sqrt{n}}$. It is true based on Lemma \ref{lem-ineq}.

Moreover, the presented functions are all m.l.b.
\end{proof}

\subsection{Proof of Theorem \ref{thm1}}\label{proof-thm1}
\begin{proof}
	Consider the candidate Lyapunov function $V(x(t)) = 2\mu(f(x(t))-f^*)$. By Lemma \ref{lem-non}, the time derivative of $V(x(t))$ can be calculated as
	\begin{align}\label{V-deri0}
		\frac{d}{dt}V(x(t)) \in \widetilde{\mathcal{L}}_{\mathcal{F}[X_u]} V(x(t)), 
	\end{align}
	For the convenience of expression, we denote $y(t)=\nabla f(x(t))$. Then, with $u(x(t))= g(\nabla f(x))$, \eqref{V-deri0} indicates that there exist $\eta_p(t) \in \mathcal{F}[g_p\circ \nabla f](x) \subset \mathcal{F}[g_p](y)$ and $\eta_q(t) \in \mathcal{F}[g_q\circ \nabla f](x) \subset \mathcal{F}[g_q](y)$ such that 
	\begin{align*}
		\frac{d}{dt}V(x(t)) &= -2\mu y^T(\eta_p+\eta_q) \\
		&\leq -2\mu \sigma \|y\|_2^{1+p}-2\mu\rho \|y\|_2^{q+1} \\
		& \leq -2\mu\sigma (V(x)))^{\frac{p+1}{2}}-2\mu\rho (V(x))^{\frac{q+1}{2}}
	\end{align*}
	where the first inequality is due to Assumption \ref{ass-gn} and the second inequality holds by virtue of \eqref{PL} in Assumption \ref{ass-f2}. Then, by Lemma \ref{lem-fixed}, any solution $x(t)$ of the nominal system \eqref{pro-u} will converge to $\mathbb{X}^*$ in fixed time bounded by \eqref{Tf-nominal}. 
\end{proof}

\subsection{Proof of Proposition \ref{prop-ge}}\label{proof-ge}
\begin{proof}
	Consider the same Lyapunov function $V(x(t))$ as that of Theorem \ref{thm1}. Denote $y(t) = \nabla f(x(t))$. For the case $g(y)=\alpha \frac{ye^{\|y\|_2}}{\|y\|_2}$, the time derivative of $V(x(t))$ is given by 
	\begin{align*}
		\frac{d}{dt}V(x(t)) &= -2\mu\alpha y^T(\frac{ye^{\|y\|_2}}{\|y\|_2}) \\
		&= -2\mu\alpha \|y\|_2 e^{\|y\|_2} \\
		&\leq -2\mu\alpha (2\mu(f(x)-f^*))^{\frac{1}{2}}e^{(2\mu(f(x)-f^*))^{\frac{1}{2}}}\\
		&= -2\mu \alpha \sqrt{V} e^{\sqrt{V}}.
	\end{align*}
	When $V(x(t))\neq 0$, we obtain 
	\begin{align*}
		\int_{V(x_0)}^{V(x(t))}\frac{1}{\sqrt{V}}e^{-\sqrt{V}}dV \leq -2\mu\alpha \int_{0}^{t}dt, 
	\end{align*}
	which implies that 
	\begin{align}\label{get}
		t \leq \frac{e^{-\sqrt{V(x(t))}}-e^{-\sqrt{V(x_0)}}}{\alpha\mu}.
	\end{align}
	Since $e^{-\sqrt{V(x(t))}}$ is lower bounded, then there exists a finite time $T_f$ such that $V(x(T_f))=0$ and we further have 
	\begin{align}\label{Tge}
		T_f \leq \frac{1- e^{-\sqrt{2\mu(f(x_0)-f^*)}}}{\alpha\mu}\leq \frac{1}{\alpha\mu}.
	\end{align}
	Hence, any solution $x(t)$ will converge to $\mathbb{X}^*$ in fixed time bounded by $\frac{1}{\alpha\mu}$. The remainder for the case $g(y)=\alpha\text{sign}(y)\odot e^{|y|}$ can be shown similarly by using inequality \eqref{ge-property2}. 
\end{proof}

\subsection{Proof of Theorem \ref{thm-gnd}}\label{proof-gnd}
\begin{proof}
	Consider the Lyapunov function $V(x(t)) = 2\mu(f(x(t))-f^*)$. By Lemma \ref{lem-non}, the time derivative of $V(x(t))$ can be calculated as
		\begin{align}\label{V-deri}
			\frac{d}{dt}V(x(t)) \in  \widetilde{\mathcal{L}}_{\mathcal{F}[X_{g,d}]} V(x(t)).
		\end{align}
		where $X_{g,d}(t,x) = g(\nabla f(x))+d(x,t)$. By the sum rule, we have 
		\begin{align*}
			\mathcal{F}[X_{g,d}](x,t) &=\mathcal{F}[g\circ \nabla f+d](x,t) \\
			& \subseteq  \mathcal{F}[g\circ \nabla f](x)+\mathcal{F}[d](x,t). 
	\end{align*}
	For the convenience of expression, we denote $y(t)=\nabla f(x(t))$. Then, with protocol \eqref{prot-gd}, \eqref{V-deri} indicates that there exist $\eta_1(t) \in \mathcal{F}[g_1\circ \nabla f](x) \subseteq \mathcal{F}[g_1](y)$, $\eta_q(t) \in \mathcal{F}[g_q\circ \nabla f](x) \subseteq \mathcal{F}[g_q](y)$ and $\eta_{d}(t) \in \mathcal{F}[d](x,t) $ such that 
	\begin{align*}
		\frac{d}{dt}V(x(t)) &= -2\mu y^T(\eta_1+ \eta_q(t) - \eta_{d}(t)) \\
		&\leq -2\mu((\sigma -\bar{d} ) \|y\|_2+ \rho \|y\|_2^{q+1} - \|y\|_2 \|x-[x]^*\|_2) \\
		&\leq -2\mu((\sigma -\bar{d} ) \|y\|_2-\rho \|y\|_2^{q+1}+ \frac{\epsilon}{2\sqrt{\mu}} \|y\|_2^2), 
	\end{align*}
	where the first inequality is due to Assumptions \ref{ass-gn}, \ref{ass-d} and the second inequality holds by virtue of \eqref{QG2}. In the next, two scenarios will be considered by investigating the value of $\|y\|_2$. When $\|y\|_1>1$, $\|y\|_2^{q+1}>\|y\|_2^2 $ and hence we obtain 
	\begin{align*}
		\frac{d}{dt}V(x(t)) \leq -2\mu (\sigma -\bar{d} ) \|y\|_2-2\mu (\rho-\frac{\epsilon}{2\sqrt{\mu}}) \|y\|_2^{q+1}.
	\end{align*}
	On the other side, when $\|y\|_1\leq 1$, $\|y\|_2\geq \|y\|_2^2 $ and we have 
	\begin{align*}
		\frac{d}{dt}V(x(t)) \leq -2\mu(\sigma -\bar{d}- \frac{\epsilon}{2\sqrt{\mu}}) \|y\|_2-2\mu\rho\|y\|_2^{q+1}.
	\end{align*}
	Let $k_1 =\sigma -\bar{d}- \frac{\epsilon}{2\sqrt{\mu}} $ and $k_2=\rho-\frac{\epsilon}{2\sqrt{\mu}}$.  For both cases, it satisfies that 
	\begin{align*}
		\frac{d}{dt}V(x(t)) &\leq -2\mu k_1 \|y\|_2-2\mu k_2 \|y\|_2^{q+1}\\
		& = -2\mu k_1 (V(x))^{\frac{1}{2}}-2\mu k_2 (V(x))^{\frac{q+1}{2}}. 
	\end{align*}
	Then, referring to Lemma \ref{lem-fixed}, if $k_1>0$ and $k_2>0$, any solution $x(t)$ of the perturbed system \eqref{pro-u} will converge to $\mathbb{X}^*$ in fixed time bounded by \eqref{Tf}. 
\end{proof}

\subsection{Proof of Theorem \ref{thm-gn2}}\label{proof-gn2}
\begin{proof}
	Consider the candidate Lyapunov function $V(x(t)) = \|\nabla f(x)\|_2^2$. By Lemma \ref{lem-non}, the time derivative of $V(x(t))$ can be calculated as
		\begin{align}\label{V-derin}
			\frac{d}{dt}V(x(t)) \in \widetilde{\mathcal{L}}_{\mathcal{F}[X_u]} V(x(t)).
	\end{align}
	Denote $y(t)=\nabla f(x(t))$. Then, with protocol \eqref{prot-gn2}, \eqref{V-derin} indicates that there exist $\eta_p(t) \in \mathcal{F}[g_p\circ \nabla f](x) \subset \mathcal{F}[g_p](y)$ and $\eta_q(t) \in \mathcal{F}[g_q\circ \nabla f](x) \subset \mathcal{F}[g_q](y)$ such that 
	\begin{align*}
		\frac{d}{dt}V(x(t)) &= -2y^T\nabla^2 f(x)(\nabla^2 f(x))^{-1}(\eta_p+\eta_q) \\
		&= -2y^T(\eta_p+\eta_q) \\
		&\leq -2 \sigma \|y\|_2^{2-p}-2\rho \|y\|_2^{q+1} \\
		& = -2 \sigma (V(x))^{1-\frac{p}{2}}-2\rho (V(x))^{\frac{q+1}{2}}
	\end{align*}
	where the first inequality is due to Assumption \ref{ass-gn}. Then, by Lemma \ref{lem-fixed}, any solution $x(t)$ of the nominal system \eqref{pro-u} will converge to $\mathbb{X}^*$ in fixed time bounded by \eqref{Tf-nominal2}. 
	
	Consider the trajectory $\mathbb{X}(x_0) = \{x(t)| t\in [0 \ T(x_0)]\}$. Since $x(t)$ is continuous over $[0 \ T(x_0)]$, $\mathbb{X}(x_0)$ is compact. As $f(x)$ is twice continuously differentiable and $\nabla^2 f(x)\succ 0$, the minimal eigenvalue denoted by $\lambda_{2}(\nabla^2 f(x))$ is continuous and attains its minimum $\underline{\lambda}>0$ on $\mathbb{X}(x_0)$. Since $(\nabla^2 f(x))^{-1}\succ 0$ and $\|(\nabla^2 f(x))^{-1}\| = \lambda_{2}^{-1}(\nabla^2 f(x))$, $\|(\nabla^2 f(x))^{-1}\| \leq \frac{1}{\underline{\lambda}}$. By Assumption \ref{ass-gn}, $g(y)$ is m.l.b. Hence, $u(x(t))$ is measurable and locally essentially bounded over $[0 \ T(x_0)]$. That is to say, the nominal dynamics \eqref{pro-u} with \eqref{prot-gn2} is well-defined. 
\end{proof}

\subsection{Proof of Theorem \ref{thm-regret}}\label{proof-regret}
\begin{proof}
	For the nominal system \eqref{pro-u}, consider $V(x(t))=f(x(t))-f^*$ with $V_0 = V(x_0)$. Denote $y(t)=\nabla f(x(t))$. Then the time derivative of $V(x(t))$ and the regret function with different protocols will be investigated. 
	1) Regret estimation for $g_0$:	With direct computation, one has that 
	\begin{align}
		\frac{d}{dt}V(x(t))= - y^Tg_1(y) = - \|y\|_2^2 \leq -(2\mu)V(x(t)),
	\end{align}
	which implies that
	\begin{align}\label{Vg0}
		V(x(t)) \leq V_0 e^{-2\mu t}.
	\end{align}
	As $V(x(t))$ will converge to zero exponentially, we can regard the convergence time $T(x_0)=+\infty$. Integrating both sides of \eqref{Vg0} from $t=0$ to $t=+\infty$ gives that 
	\begin{align}
		\mathcal{R}(T(x_0),x_0) \leq \frac{V_0}{2\mu}.
	\end{align}
	
	2) Regret estimation for $g_p$: From the proof of Proposition \ref{prop-FT}, we obtain
	\begin{align}\label{Vgp}
		V(x(t)) \leq (V_0^{1-\alpha}-a(1-\alpha) t)^{\frac{1}{1-\alpha}}
	\end{align}
	with $\alpha= \frac{p+1}{2}$ and $a=(2\mu)^{\alpha}$ for any $t\geq0$. Integrating both sides of \eqref{Vgp} from $t=0$ to $t=T(x_0)$ referring to \eqref{Tf2-gn} gives that 
	\begin{align}\label{Rgp}
		\mathcal{R}(T(x_0),x_0) \leq\frac{V_0^{2-\alpha}}{a(2-\alpha)}.
	\end{align}
	
	3) Regret estimation for $g_{p,q}$: Following the proof of Theorem \ref{thm1}, it can be derived that
	\begin{align*}
		\frac{d}{dt}V(x(t)) & \leq -(2\mu)^{\alpha}V^{\alpha}(x) -(2\mu)^{\beta}V^{\beta}(x) \\
		& \leq -a V^{\alpha}(x) -b V^{\beta}(x)
	\end{align*}
	with $\alpha= \frac{p+1}{2}$ and $\beta=\frac{q+1}{2}$ for any $t\geq0$. Following the previous part 2), \eqref{Rgp} holds when $V_0\leq 1$. When $V_0> 1$, motivated by the proof of \cite[Thm. 3]{Budhraja2021}, one can consider two stages. One is that $x(t)$ converges to $S=\{x: V(x)\leq 1 \}$ in fixed time $T_1 \leq \frac{1}{b(\beta-1)}$ dominated by the dynamics $\dot{V}(x(t))\leq -b V^{\beta}(x)$ and it holds that
	\begin{align}\label{Vgpq}
		V(x(t)) \leq \frac{V_0}{(1+bV_0^{\beta-1}(\beta-1)t)^{\frac{1}{\beta-1}}}.
	\end{align}
	Then, \eqref{Rgp} can be applied when $V(x(t)) \leq 1$. When $V_0>1$ and $\beta\neq 2$, combining \eqref{Vgpq} with \eqref{Rgp}, it can be derived that
	\begin{align}\label{Rfxt}
		\mathcal{R}(T(x_0),x_0) \leq \frac{1}{a(2-\alpha)} + \frac{(1+V_0^{\beta-1})^{\frac{\beta-2}{\beta-1}}-1}{bV_0^{\beta-2}(\beta-2) }.
	\end{align}
	Similarly, when $V_0>1$ and $\beta= 2$, one can obtain 
	\begin{align}\label{Rfxt-2}
		\mathcal{R}(T(x_0),x_0) \leq \frac{1}{a(2-\alpha)} + \frac{\text{ln}(1+V_0)}{b }.
	\end{align}
	
	Define
	\begin{align}\label{Fdef}
		F(V_0, \beta) \triangleq \frac{(1+V_0^{\beta-1})^{\frac{\beta-2}{\beta-1}}-1}{V_0^{\beta-2}(\beta-2) }.
	\end{align} 
When $\beta>2$, $F(V_0, \beta)$ can be relaxed by 
	\begin{align}
		F(V_0, \beta) < \frac{1}{\beta-2}
	\end{align}
	using inequality $(1+z)^{c} <1+z^c $ for any $z>0$ and $c<1$. When $\beta \in (1,2)$, one can estimate $F(V_0, \beta)$ as follows
	\begin{align}
		F(V_0, \beta) = \frac{1}{2-\beta}(V_0^{2-\beta}-\frac{1}{(\frac{1}{V_0^{\beta-1}}+1)^{\frac{2-\beta}{\beta-1}}})
		< \frac{V_0^{2-\beta}}{2-\beta}.
	\end{align}
	
	4) Regret estimation for $g_{e}$: From \eqref{get}, it can be derived that 
	\begin{align}\label{Vge}
		V(x(t)) \leq \frac{1}{2\mu} \text{ln}^2(\mu t + e^{-\sqrt{2\mu V_0}}).
	\end{align}
	Let $w = e^{-\sqrt{2\mu V_0}}$ and $y= \ln (\mu t + w)$. Integrating both sides of \eqref{Vge} from $t=0$ to $t=T(x_0)$ with \eqref{Tge} yields
	\begin{align*}
		\mathcal{R}(T,x_0) &\leq \frac{1}{2\mu} \int_{0}^{T}\text{ln}^2(\mu t + w) dt\\
		& =\frac{1}{2\mu^2} \int_{\ln (w)}^{\ln (\mu T+w)} y^2 de^y\\
		& =\frac{1}{2\mu^2} ((y-1)^2+1)e^y\bigg|_{\ln (w)}^{\ln (\mu T+w)} \\
		& =\frac{1}{2\mu^2} ((\ln s-1)^2+1) s\bigg|_{w}^{\mu T+w}  \\
	\end{align*}
	Let $h(s)=((\ln s-1)^2+1) s$, which satisfies $\nabla h(s) = \ln^2s \geq 0$. Since $\mu T(x_0) +w \leq 1$, we further give 
	\begin{align*}
		\mathcal{R}(T(x_0),x_0) &\leq  \frac{h(1)-h(w)}{2\mu^2}\\
		&=\frac{2-((1+\sqrt{2\mu V_0})^2+1 )e^{-\sqrt{2\mu V_0}}}{2\mu^2} \leq \frac{1}{\mu^2}.
	\end{align*}
\end{proof}
\bibliographystyle{ieeetr}

\begin{thebibliography}{10}



%
%
%
%
\bibitem{Wang2011}
J. Wang and N. Elia, “A control perspective for centralized and distributed
convex optimization,” in {\em Proc. IEEE Conf. Decision Control and Eur. Control Conf.}, Orlando, FL, Dec. 2011, pp. 3800-3805.

\bibitem{Liu2017}
Q. Liu, S. Yang, and J. Wang, “A collective neurodynamic approach to distributed constrained optimization,” {\em IEEE Trans. Neural Netw.}, vol. 28, no. 8, pp. 1747-1758, 2017.


\bibitem{cortes2006}
J.~Cort\'{e}s, “Finite-time convergent gradient flows with applications to network consensus,” {\em Automatica}, vol.~42, no.~11, pp. 1993-2000, 2006.

\bibitem{Xu1998} 
C. Xu and J. Prince, “Snakes, shapes, and gradient vector flow,” {\em IEEE Trans. Image Process.}, vol. 7, no. 3, pp. 359-369, Mar. 1998.

\bibitem{Su2014}
W. Su, S. Boyd, and E. Candes, “A differential equation for modeling nesterov’s accelerated gradient method: Theory and insights,” in {\em Proc. Adv. Neural Inf. Process. Syst.}, pp. 2510-2518, 2014.

\bibitem{Wibisono2016}
A. Wibisono, A. C. Wilson, and M. I. Jordan, “A variational perspective on accelerated methods in optimization,” {\em Nat. Acad. Sci.}, vol. 113, no. 47, pp. E7351-E7358, 2016.

\bibitem{Attouch2018}
H. Attouch, Z. Chbani, J. Peypouquet, and P. Redont, “Fast convergence of inertial dynamics and algorithms with asymptotic vanishing viscosity,” {\em Math. Program.}, vol. 168, no. 1-2, pp. 123-175, 2018.

\bibitem{Vassilis2018}
A. Vassilis, A. Jean-Fran\c cois, and D. Charles, “The differential inclusion modeling FISTA algorithm and optimality of convergence rate in the case $b \leq 3$,” {\em SIAM J. Optim.}, vol. 28, pp. 551-574, 2018. 

\bibitem{Sebbouh2020}
O. Sebbouh, C. Dossal, and A. Rondepierre, “Convergence rates of damped inertial dynamics under geometric conditions and perturbations,” {\em SIAM J. Optim.}, vol. 30, pp. 1850-1877, 2020.

\bibitem{Serna2021}
J. M. Sanz Serna, and K. C. Zygalakis, “The connections between Lyapunov functions for some optimization algorithms and differential equations,” {\em SIAM. J. Numer. Anal.}, vol. 59, no. 3, pp. 1542-1565, 2021.

\bibitem{Bhat2000}
S. P. Bhat and D. S. Bernstein, “Finite-time stability of continuous autonomous systems,”
{\em SIAM J. Control Optim.}, vol. 38, no. 3, pp. 751-766, 2000.

\bibitem{Yu2005}
S. Yu, X. Yu, B. Shirinzadeh, and Z. Man, “Continuous finite-time control for robotic manipulators with terminal sliding mode,” {\em Automatica}, vol. 41, no. 11, pp. 1957-1964, 2005.
\bibitem{Shen2012}
Y. Shen and Y. Huang, “Global finite-time stabilisation for a class of nonlinear systems,” {\em
Int. J. Syst. Sci.}, vol. 43, no. 1, pp. 73-78, 2012.
\bibitem{Aouiti2021FT}
C. Aouiti and M. Bessifi, “Periodically intermittent control for finite-time synchronization of delayed quaternion-valued neural networks,” {\em Neural. Comput. Appl.}, vol. 33, pp. 6527-6547, 2021.



\bibitem{Romero2020}
O. Romero and M. Benosman, “Finite-time convergence in continuous-time optimization,” in {\em International Conference on Machine Learning}, PMLR, 2020, pp. 8200-8209.

\bibitem{Chen2020sign}
F. Chen and W. Ren, “Sign projected gradient flow: A continuous time approach to convex optimization with linear equality constraints,” {\em Automatica}, vol. 120, p. 109156, 2020.




\bibitem{Wei2022}
Y. Wei, Y. Chen, X. Zhao, and J. Cao, “Analysis and synthesis of gradient algorithms based on fractional-order system theory,” {\em IEEE Trans. Syst., Man, Cybern., Syst.}, vol. 53, no. 3, pp. 1895-1906, March 2023.

\bibitem{Zhou2019tac}
J. Zhou, X. Wang, S. Mou, and B. D. Anderson, “Finite-time distributed linear equation solver for solutions with minimum $l_1$-norm,” {\em IEEE Trans. Autom. Control}, vol. 65, no. 4, pp. 1691-1696, 2020.



\bibitem{Shi2022cyber}
X. Shi, X. Xu, X. Yu, and J. Cao, “Finite-time convergent primal-dual gradient dynamics with applications to distributed optimization,” {\em IEEE Trans. Cybern.}, vol. 53, no. 5, pp. 3240-3252, 2023.

\bibitem{Shi2022TAC}
X. Shi, G. Wen, J. Cao, and X. Yu, “Finite-time distributed average tracking for multi-agent optimization with bounded inputs,” {\em IEEE Trans. Autom. Control}, vol. 68, no. 8, pp. 4948-4955, 2023.


\bibitem{Shi2023}
X. Shi, G. Wen and X. Yu, “Finite-time convergent algorithms for time-varying distributed optimization,” {\em IEEE Control Syst. Lett.}, vol. 7, pp. 3223-3228, 2023.


\bibitem{Polyakov2012}
A. Polyakov, “Nonlinear feedback design for fixed-time stabilization of linear control systems,” {\em IEEE Trans. Autom. Control}, vol. 57, pp. 2106-2110, 2012.

\bibitem{LiuCAA2022}
Y. Liu, H. Li, Z. Zuo, X. Li, and R. Lu, “An overview of finite/fixed-time control and its application in engineering systems,” {\em IEEE/CAA J. Autom. Sinica}, vol. 9, no. 12, pp. 1–15, 2022.

\bibitem{Aouiti2019}
C. Aouiti, E. A. Assali, and Y. E. Foutayeni, “Finite-time and fixed-time synchronization of
inertial Cohen–Grossberg-type neural networks with time varying delays,” {\em Neural Processing Letters}, vol. 50, pp. 2407-2436, 2019.
\bibitem{Alimi2019}
A. M. Alimi, C. Aouiti, and E. A. Assali, “Finite-time and fixed-time synchronization of a class of inertial neural networks with multi-proportional delays and its application to secure communication,” {\em Neurocomputing}, vol. 332, pp. 29-43, 2019.
\bibitem{Cao2017}
J. Cao and R. Li, “Fixed-time synchronization of delayed memristor-based recurrent neural networks,” {\em Sci. China Inf. Sci.}, vol. 60, no. 3, 2017.
\bibitem{Aouiti2020}
C. Aouiti and F. Miaadi, “A new fixed-time stabilization approach for neural networks with time-varying delays,” {\em Neural. Comput. Appl.}, vol. 32, pp. 3295-3309, 2020.
\bibitem{Aouiti2020CSSP}
C. Aouiti, M. Bessifi, and X. Li, “Finite-time and fixed-time synchronization of complex-valued recurrent neural networks with discontinuous activations and time-varying delays,” {\em Circuits, Syst., Signal Process.}, vol. 39,
no. 11, pp. 5406–5428, 2020.
\bibitem{Aouiti2021}
C. Aouiti, Q. Hui, H. Jallouli, and E. Moulay, “Sliding mode control-based fixed-time stabilization and synchronization of inertial neural networks with time-varying delays,” {\em Neural. Comput. Appl.}, vol. 33, pp. 11555-11572, 2021.



\bibitem{Garg2021tac}
K. Garg and D. Panagou, “Fixed-time stable gradient flows: Applications to continuous-time optimization,” {\em IEEE Trans. Autom. Control}, vol. 66, no. 5, pp. 2002-2015, 2020.

\bibitem{Guo2022CAA}
L. Guo, X. Shi, and J. Cao, “Exponential convergence of primal-dual dynamical system for linear constrained
optimization,” {\em IEEE/CAA J. Autom. Sinica}, vol. 9, no. 4, pp. 745–748, 2022.

\bibitem{Budhraja2021}
P. Budhraja, M. Baranwal, K. Garg, and A. Hota, “Breaking the convergence barrier: Optimization via fixed-time convergent flows,” In Proceedings of the AAAI Conference on Artificial Intelligence, vol. 36, no. 6, 2022.


\bibitem{Garg2019MVI}
K. Garg, M. Baranwal, R. Gupta, and M. Benosman, “Fixed-time stable proximal dynamical system for solving MVIPs,” {\em IEEE Trans. Autom. Control}, vol. 68, no. 8, pp. 5029-5036, 2023.


\bibitem{Ju2021Tcyber}
X. Ju, D. Hu, C. Li, X. He, and G. Feng, “A novel fixed-time converging neurodynamic approach to mixed variational inequalities and applications,” {\em IEEE Trans. Cybern.}, vol. 52, no. 12, pp. 12942-12953, 2022.



\bibitem{HeTNNLS2021}
X. He, H. Wen, and T. Huang, “A fixed-time projection neural network for solving $L_1$-minimization problem,” {\em IEEE Trans. Neural Netw.}, vol. 33, no. 12, pp. 7818-7828, Dec. 2022.


\bibitem{Wu2022SMC}
Z. Wu, Z. Li and J. Yu, “Designing zero-gradient-sum protocols for finite-time distributed optimization problem,” {\em IEEE Trans. Syst., Man, Cybern., Syst.}, vol. 52, no. 7, pp. 4569-4577, July 2022.







\bibitem{Shi2020auto}
X. Shi, X. Yu, J. Cao, and G. Wen, “Continuous distributed algorithms for solving linear equations in finite time,” {\em Automatica}, vol. 113, p. 108755, 2020.

\bibitem{Karimi2016}
H. Karimi, J. Nutini, and M. Schmidt, “Linear convergence of gradient and proximal-gradient methods under the Polyak-Łojasiewicz condition,” in {\em Proc. Eur. Conf. Mach. Learn.}, Sep. 2016, pp. 795–811. 

\bibitem{Yi2022CAA}
X. Yi, S. Zhang, T. Yang, T. Chai, and  K. H. Johansson,  “A primal-dual SGD algorithm for distributed nonconvex optimization,” {\em IEEE/CAA J. Autom. Sinica}, vol. 9, no. 5, pp. 812–833, 2022. 


 
%
 
 \bibitem{Shi2020TAC}
 X. Shi, J. Cao, X. Yu, and G. Wen, “Finite-time stability for network systems with discontinuous dynamics over signed digraphs,” {\em IEEE Trans. Autom. Control}, vol. 65, no. 11, pp. 4874-4881, 2020. 
 
 \bibitem{Xiao2009}
 F. Xiao, L. Wang, J. Chen, and Y. Gao, ``Finite-time formation control for multi-agent systems,'' {\em Automatica}, vol. 45, no. 11, pp. 2605-2611, 2009.
 
 \bibitem{Shi2019TNSE}
 X. Shi, J. Cao, G. Wen, and X. Yu, “Finite-time stability for network systems with nonlinear protocols over signed digraphs,” {\em IEEE Trans. Netw. Sci. Eng.}, vol. 7, no. 3, pp. 1557-1569, 2020.
 
 \bibitem{Zuo2014}
 Z. Zuo and L. Tie, ``A new class of finite-time nonlinear consensus protocols for multi-agent systems,'' {\em Int. J. Control}, vol. 87, no. 2, pp. 363-370, 2014.
 
 \bibitem{Cortes2008}
 J.~Cort\'{e}s, “Discontinuous dynamical systems: A tutorial on solutions,
 nonsmooth analysis, and stability,” {\em IEEE Control Syst. Mag.}, vol. 28, no. 3, pp. 36-73, 2008.
 
  \bibitem{Yu2021} 
 X. Yu, Y. Feng and Z. Man, “Terminal sliding mode control – An overview,” {\em IEEE Open J. Ind. Electron. Soc.}, vol. 2, pp. 36-52, 2021.
 
 
\bibitem{Parikh2013}
 N. Parikh and S. Boyd, “Proximal algorithms,” {\em Found. Trends Optim.}, vol. 1, no. 3, pp. 123–231, 2014.
 

 \bibitem{Themelis2018}
A. Themelis, L. Stella, and P. Patrinos, “Forward-backward envelope for the sum of two nonconvex functions: Further properties and nonmonotone line-search algorithms,” {\em SIAM J. Optim.}, vol. 28, no. 3, pp. 2274–2303, 2018.

\bibitem{Hassan2021AUTO}
S. Hassan-Moghaddam and M. R. Jovanovi\'{c}, “Proximal gradient flow and Douglas–Rachford splitting dynamics: Global exponential stability via integral quadratic constraints,” {\em Automatica}, vol. 123, p. 109311, 2021. 

 
\bibitem{Goldberg1987}
M. Goldberg, “Equivalence constants for $l_p$ norms of matrices,” {\em Lin. Multilin. Algebra}, vol. 21, no. 2, pp. 173-179, 1987.



\bibitem{Wang2010}
L.~Wang and F.~Xiao, ``Finite-time consensus problems for networks of dynamic agents,'' {\em IEEE Trans. Autom. Control}, vol.~55, no.~4, pp.~950-955, 2010.

\bibitem{Hui2008}
Q. Hui, W. Haddad, and S. Bhat, ``Finite-time semistability and consensus for nonlinear dynamical networks,'' {\em IEEE Trans. Autom. Control}, vol. 53, no. 8, pp. 1887-1890, 2008.
 
 

\bibitem{Yang2015}
M. Yang and C. Y. Tang, ``A distributed algorithm for solving general linear equations over networks,'' {\em Proc. IEEE Conf. Decis. Control}, pp. 3580-3585, Dec. 2015.

 
\bibitem{Wand2019arc}
 P. Wang, S. Mou, J. Lian, and W. Ren, “Solving a system of linear equations: From centralized to distributed algorithms,” {\em Annu. Rev. Control}, vol. 47, pp. 306-322, 2019.
 

 
\bibitem{Lee2016}
S. Lee, A. Ribeiro and M. M. Zavlanos, “Distributed continuous-time online optimization using saddle-point methods,” in {\em Proc. 55th IEEE Conf. Decis. Control}, Dec. 2016, pp. 4314-4319.

\bibitem{Wood2013}
A. J. Wood, B. F. Wollenberg, and G. B. Sheble, “Power generation, operation, and control,” New York, NY: Wiley, 2013.

\bibitem{Cherukuri2015}
A. Cherukuri and J. Cort\'es, “Distributed generator coordination for initialization and anytime optimization in economic dispatch,” {\em IEEE Trans. Control Netw. Syst.}, vol. 2, no. 3, pp. 226-237, 2015.



























%







%
%















%




\end{thebibliography}

\end{document}